\documentclass[pdflatex,sn-mathphys]{sn-jnl}

\usepackage{tikz}
\usetikzlibrary{arrows}

\jyear{2021}%

\theoremstyle{thmstyleone}%
\newtheorem{theorem}{Theorem}
\newtheorem{proposition}[theorem]{Proposition}%
\newtheorem{lemma}[theorem]{Lemma}

\newtheorem{corollary}[theorem]{Corollary}

\theoremstyle{thmstyletwo}%
\newtheorem{example}{Example}%
\newtheorem{remark}{Remark}%
\newtheorem*{note}{Note}%
\newtheorem{question}{Question}
\newtheorem*{redi_proof}{Proof}%

\theoremstyle{remark}
\newtheorem{construction}{Construction}[section]

\theoremstyle{thmstylethree}%
\newtheorem{definition}{Definition}%

\begin{document}

\title[A simplicial category for higher correspondences]{A simplicial category for higher correspondences}


\author{\fnm{Redi} \sur{Haderi}}\email{redi.haderi@bilkent.edu.tr}

\affil{\orgdiv{Department of Mathematics}, \orgname{Bilkent University}, \orgaddress{\street{Cankaya}, \city{Ankara}, \postcode{06800}, , \country{Turkey}}}


\abstract{In this work we propose a realization of Lurie's prediction that inner fibrations $p: X \rightarrow A$ are classified by $A$-indexed diagrams in a ``higher category" whose objects are $\infty$-categories, morphisms are correspondences between them and higher morphisms are higher correspondences. We will obtain this as a corollary of a more general result which classifies all simplicial maps between ordinary simplicial sets in a similar fashion.
    
    Correspondences between simplicial sets (and $\infty$-categories) are a generalization of the concept of profunctor (or bimodule) pertaining to categories. While categories, functors and profunctors are organized in a double category, we will exhibit simplicial sets, simplicial maps, and correspondences as part of a simplicial category. This allows us to make precise statements and provide proofs. Our main tool is the  language of double categories, which we use in the context of simplicial categories as well. }

\keywords{simplicial category, correspondence, double category, double colimit}



\maketitle

\section*{Introduction and summary}

The main object of study in this work is the notion of correspondence between simplicial sets and its role in classifying maps of simplicial sets. We will do so using tools from double category theory developed in \cite{grandis1999limits}. While the subject is certainly of independent interest, and perhaps fundamental in nature, its study can be strongly motivated from $\infty$-category\footnote{We use this term according to Lurie's book \cite{lurie2009higher}.} theory as well: correspondences between $\infty$-categories are expected to fit in a higher structure which ought to serve as a classifier for inner fibrations. 

In order to precisely state the ideas in the previous paragraph, we need to revisit certain aspects of a classical subject: correspondences between categories, also known as profunctors. Recall that for two categories $\mathcal{C}$ and $\mathcal{D}$ a $(\mathcal{C}, \mathcal{D})$-profunctor is a functor 
$$u : \mathcal{C}^{op} \times \mathcal{D} \rightarrow \textbf{Set}$$
into the category of sets. 

Equivalently (as observed by B\'{e}nabou, see \cite{street2001powerful}), a profunctor $u$ can be presented via its collage $\textbf{col}(u)$, which is a category whose objects are those of $\mathcal{C}$ and $\mathcal{D}$ and sets of morphisms defined as 
$$\textbf{col}(u)(a,b) = \begin{cases}
\mathcal{C}(a,b), & \text{if} \ \  a,b \in \mathcal{C} \\
\mathcal{D}(a,b), & \text{if} \ \  a,b \in \mathcal{D} \\
u(a,b), & \text{if} \ \  a \in \mathcal{C}, b \in \mathcal{D} \\
\emptyset, & \text{if} \ \  a \in \mathcal{D}, b \in \mathcal{C} \\
\end{cases}$$
$\textbf{col}(u)$ is naturally equipped with a map $p : \textbf{col}(u) \rightarrow \Delta^1$, where $\Delta^1 = \{ 0 \rightarrow 1 \}$ is the interval category, such that $p^{-1}(0) \cong \mathcal{C}$ and $p^{-1}(1) \cong \mathcal{D}$. Conversely, any category $\mathcal{U}$ equipped with such a functor $p$ is the collage of a profunctor. 

Profunctors can be tensored (a.k.a. composed) with the usual coend formula.
It is common to organize categories, profunctors and transformations between them in a 2-category $\textbf{Cat}_{prof}$. However, as pointed out in \cite{grandis1999limits}, it is more advantageous, and arguably natural, to form a double category \textbf{Prof} with profunctors as horizontal morphisms, functors as vertical morphisms and 2-cells
\begin{center}
        \begin{tikzpicture}
\node (v1) at (-2.5,3) {$\mathcal{C}_0$};
\node (v4) at (-2.5,1) {$\mathcal{D}_0$};
\node (v2) at (-0.5,3) {$\mathcal{C}_1$};
\node (v3) at (-0.5,1) {$\mathcal{D}_1$};
\draw[->]  (v1) edge node[above, font = \scriptsize] {$u$} (v2);
\draw[->]  (v2) edge node[right, font = \scriptsize] {$G$} (v3);
\draw[->]  (v1) edge node[left, font = \scriptsize] {$F$} (v4);
\draw[->]  (v4) edge node[below, font = \scriptsize] {$v$}  (v3);
\node (v5) at (-1.5,2.6) {};
\node (v6) at (-1.5,1.4) {};
\draw[->]  (v5) edge[double] node[right] {$\alpha$} (v6);
        \end{tikzpicture}
    \end{center}
given by morphisms of profunctors $\alpha$ which respect $F$ and $G$. 

This double categorical formalism is particularly relevant in virtue of the following. Let $F : \mathcal{X} \rightarrow \mathcal{A}$ be a functor and for objects $a \in \mathcal{A}$ let $\mathcal{X}_a = F^{-1}(a)$ be the fibres. In light of the above, for a morphism $f: a \rightarrow b$ the fibre $\mathcal{X}_f = F^{-1}(f)$ presents an $(\mathcal{X}_a, \mathcal{X}_b)$-profunctor. This way the functor $F$ can be, borrowing a term from Lurie (\cite{lurie2009higher}), straightened into a lax 2-functor 
$$\chi_F : \mathcal{A} \rightarrow \textbf{Cat}_{prof}$$

It is expected of $F \mapsto \chi_F$ to be functorial, i.e. there should be a functor 
$$\chi : \textbf{Cat}/\mathcal{A} \rightarrow [\mathcal{A}, \textbf{Cat}_{prof}]_{lax}$$
from the slice category over $\mathcal{A}$ into a category whose objects are lax diagrams of profunctors. But the 2-categorical notion of transformation between lax diagrams is not appropriate for $\chi$ to be a functor. Rather, the notion of vertical transformation coming from double category theory is the correct choice of morphism. 

Moreover, $\chi$ turns out to be an equivalence of categories. We observe that its weak inverse, which is provided by an unstraightening process, is naturally understood in double category theory as a double colimit (in the sense of \cite{grandis1999limits}). Hence our first result is a double categorical restatement of the known  fact (\cite{street2001powerful}) that functors are classified by lax diagrams of profunctors.

\begin{theorem}
Let $\mathcal{A}$ be a category. The pair $(\chi, \textbf{dcolim})$
$$\chi : \textbf{Cat}/\mathcal{A} \leftrightarrows [\mathcal{A}, \textbf{Prof}_h]_{lax} : \textbf{dcolim}$$
 is an equivalence of categories between the slice category of functors over $\mathcal{A}$ and the category of lax $\mathcal{A}$-indexed diagrams of profunctors with vertical transformations as morphisms. 

\end{theorem}

The main question we consider is the following.

\begin{question}\label{qprof}
 Is there a double category of simplicial sets which is analogous to \textbf{Prof}? 
\end{question}

It is easy enough to define correspondences between simplicial sets: they are maps to $\Delta^1$. However, there is no natural way to compose them, and hence we cannot form a double category of simplicial sets whose vertical morphisms are simplicial maps and horizontal morphisms are correspondences.

Nonetheless, we claim that if we expand the scope of Question \ref{qprof} to ask for a 2-fold structure rather than a double category, we can develop a meaningful double category theory for simplicial sets. This is based on the observation that even without a composition correspondences are still part of a simplicial structure.
For this we have to invoke higher correspondences. 

Let an $n$-correspondence be defined to be a map of simplicial sets $p: X \rightarrow \Delta^n$, where $\Delta^n$ is the standard $n$-simplex. Let $X_i = p^{-1}(i)$ be the fibres of $p$ at the vertices of $\Delta^n$. The structure map $p$ allows us to think of $X$ as a collage, formed from the $X_i$'s by adding new simplices in the increasing direction. In particular, $X$ itself looks and feels like an $n$-simplex. 

This idea can be formalized by constructing a simplicial category $\textbf{sSet}^\sharp$, i.e. a mapping
$$\textbf{sSet}^\sharp : \Delta^{op} \rightarrow \textbf{Cat}$$
with category of $n$-simplices the category of $n$-correspondences, i.e. the slice category
$$\textbf{sSet}^\sharp_n = \textbf{sSet}/\Delta^n$$
$\textbf{sSet}^\sharp$ turns out to be analogous to \textbf{Prof} from many different angles.

First, we may regard a simplicial category $\mathbb{E}$ as a 2-fold structure with objects of $\mathbb{E}_n$ as simplices in the horizontal direction and morphisms in $\mathbb{E}_n$ in the vertical direction.

Second, in general simplicial categories are 2-fold analogs of simplicially enriched categories (\textbf{sSet}-categories for short). By forgetting the horizontal part of a double category, meaning discarding all non-identity horizontal morphisms, we are left with a 2-category. When we forget the horizontal part of a simplicial category, i.e. by discarding all simplices  except those which are totally degenerate (meaning every face is degenerate), we are left with an \textbf{sSet}-category. The vertical part of $\textbf{sSet}^\sharp$ is precisely the usual \textbf{sSet}-category of simplicial sets. 

Lastly, each simplicial map $f: X \rightarrow A$ can be straightened into a map of simplicial categories
$$\chi_f : A \rightarrow \textbf{sSet}^\sharp$$
This mapping assigns to a simplex $\sigma \in A_n$ the $n$-correspondence $X_\sigma = f^{-1}(\sigma)$. 

Moreover, with degeneracies playing the role of identities and simplices the role of composites it is straightforward to develop a theory of double colimits for simplicial categories. We prove the following theorem in \ref{class_results}, which says that maps of simplicial sets are classified by diagrams of correspondences.

\begin{theorem}
Let $A$ be a simplicial set. The pair of functors $(\chi, \textbf{dcolim})$ 
$$\chi: \textbf{sSet}/A \leftrightarrows [A, \textbf{sSet}^\sharp] : \textbf{dcolim}$$
is an equivalence of categories.
\end{theorem}
We also briefly sketch an enhancement of this result to an equivalence of simplicially enriched categories.

As an application of the above we make precise Lurie's prediction that ``inner fibrations are classified by diagrams into a higher category whose objects are $\infty$-categories and 1-morphisms are correspondences"(see Section 2.3.1 in \cite{lurie2009higher}). We can form a simplicial category $\infty\textbf{Cat}^\sharp$ of $\infty$-categories, exactly as $\textbf{sSet}^\sharp$, as a model for this higher category. The stability properties of inner fibrations yield that a map $f: X \rightarrow A$ is an inner fibration if and only if each fiber $X_\sigma$, $\sigma \in A$, is an $\infty$-category. In light of Theorem \ref{class_simp}, $f$ is an inner fibration if and only if the classifying map $\chi_f$ factors through the inclusion $\infty\textbf{Cat}^\sharp \subset \textbf{sSet}^\sharp$. In particular there is an equivalence of categories
$$(\textbf{sSet}/A)_{inner} \cong [A, \infty\textbf{Cat}^\sharp]$$
between the full subcategory of the slice category over $A$ consisting of inner fibrations and the category of $A$-indexed diagrams of correspondences of $\infty$-categories.

\section*{Organization}

In Section \ref{intro} we introduce our key players: profunctors (correspondences), double categories and higher correspondences. In Section \ref{pattern} we present from a double categorical lens how functors are classified by  lax diagrams of profunctors. Then in Section \ref{higher_corrs} we construct a simplicial category of higher correspondences and state and prove the corresponding results.

\section*{Prerequisites, notation and conventions}

We expect the reader to be familiar with basic category theory. Some familiarity with 2-categories and simplicial sets is also assumed, although deep knowledge of simplicial homotopy theory is not needed to read this document. Our references for these subjects are the first chapter of \cite{leinster2004higher} and the introductory \cite{riehl2011leisurely} respectively. 

We will denote categories in general by calligraphic letters like $\mathcal{C}$, $\mathcal{D}$ etc., while we will name specific categories after their objects. For example \textbf{Set} will denote the category of sets, \textbf{Cat} the category of categories, \textbf{sSet} the category of simplicial sets etc. If $c$ is an object of $\mathcal{C}$ we will write $c \in \mathcal{C}$ and write $1_c$ for its identity morphism. If $c,d \in \mathcal{C}$ we denote the set of morphisms between them $\mathcal{C}(a,b)$. $\mathcal{C}^{op}$ is the opposite category. General categories $\mathcal{C}$ are always assumed to be small, while we allow large locally small categories to be part of \textbf{Cat}.

We denote by $\Delta$ the category of finite ordinals with order preserving maps as morphisms. The ordinal $[n] \in \Delta$ is simply $\{0 < 1 < \dots < n\}$. This way a simplicial set is a functor $X : \Delta^{op} \rightarrow \textbf{Set}$. The image of $[n]$ under $X$, or the set of $n$-simplices of $X$, is denoted by $X_n$. For a morphism $\theta : [m] \rightarrow [n]$ in $\Delta$ we denote by $\sigma \theta$ the action of $\theta$ on a simplex $\sigma \in X_n$.
Simplicial sets form a category with natural transformations as morphisms. 

As usual in the literature, we denote by $\Delta^n$ the simplicial set represented by $[n]$. Also, we will write $d^i : [n] \rightarrow [n+1]$ for the map defined as
$$d^i(j) = \begin{cases}
j & ,\text{if} \ \ j < i \\
j+1 & , \text{if} \ \ j > i
\end{cases}$$
and $s^i : [n +1] \rightarrow [n]$ for the map defined as
$$s^i(j) = \begin{cases}
j & ,\text{if} \ \ j \leq i \\
j-1 & , \text{if} \ \ j > i
\end{cases}$$
whenever these expressions make sense. For a simplicial set $X$, we write the images of these maps $d_i : X_{n+1} \rightarrow X_n$ and $s_i : X_n \rightarrow X_{n+1}$, and call them face and degeneracy maps respectively. 

Given the scope of our paper, the term ``simplicial category" means simplicial object in \textbf{Cat}, i.e. a map $\Delta^{op} \rightarrow \textbf{Cat}$, while we refer to simplicially enriched categories as \textbf{sSet}-categories. 

\section{Preliminaries} \label{intro}

\subsection{Correspondences} \label{correspondences}

\subsubsection{Same object, many faces} \label{profunctors}

Besides functors, another interesting notion of morphism between categories is that of a profunctor. 

\begin{definition}
A profunctor $u$ between two categories $\mathcal{C}$ and $\mathcal{D}$, or $(\mathcal{C}, \mathcal{D})$-profunctor, is a functor
$$u: \mathcal{C}^{op} \times \mathcal{D} \rightarrow \textbf{Set}$$
where \textbf{Set} is the category of sets and functions.
\end{definition}

A profunctor may be seen to record a right action of $\mathcal{C}$ and a left action of $\mathcal{D}$ which commute. If $\mathcal{C}$ and $\mathcal{D}$ are groups, i.e. single object categories with all morphisms invertible,  profunctors are also known as \textit{bisets}. If we enrich the above definition over abelian groups and let $\mathcal{C}$ and $\mathcal{D}$ be rings (meaning they have one object) then profunctors are simply \textit{bimodules}. If $\mathcal{C}$ and $\mathcal{D}$ are sets then all relations are profunctors (having as values only the empty set or the singleton set). This is why they are referred to as \textit{relators} sometimes. 

Examples are abundant. By taking $\mathcal{C}$ or $\mathcal{D}$ to be the terminal category (with one object and only the identity morphism) any presheaf or copresheaf is a profunctor. Also, any functor $F: \mathcal{C} \rightarrow \mathcal{D}$ produces a $(\mathcal{C}, \mathcal{D})$-profunctor $F^*$ given by 
$$F^*(c,d) = \mathcal{D}(Fc, d)$$
and a $(\mathcal{D}, \mathcal{C})$-profunctor $F_*$ given by
$$F_*(d,c) = \mathcal{D}(d, Fc)$$
for pairs $c \in \mathcal{C}$, $d \in \mathcal{D}$. $F^*$ may be thought of as the bimodule induced by $F$ and dually $F_*$ as the bimodule coinduced by $F$.

The perspective of interest from our point of view is one of a more combinatorial flavour. Given a $(\mathcal{C}, \mathcal{D})$-profunctor $u$ we may record all its information by forming a category $\textbf{col}(u)$ called the \textit{collage} of $u$. $\textbf{col}(u)$ has as objects (the disjoint union of) those of $\mathcal{C}$ and $\mathcal{D}$. For two objects $a, b \in \textbf{col}(u)$ the set of morphisms between them is defined to be
$$\textbf{col}(u)(a,b) = \begin{cases}
\mathcal{C}(a,b), & \text{if} \ \  a,b \in \mathcal{C} \\
\mathcal{D}(a,b), & \text{if} \ \  a,b \in \mathcal{D} \\
u(a,b), & \text{if} \ \  a \in \mathcal{C}, b \in \mathcal{D} \\
\emptyset, & \text{if} \ \  a \in \mathcal{D}, b \in \mathcal{C} \\
\end{cases}$$

In other words, the collage of a profunctor $u$ is produced by depicting the elements of the sets $u(c,d)$ for $c \in \mathcal{C}$, $d \in \mathcal{D}$ as actual arrows $c \rightarrow d$, on top of the arrows already existing in $\mathcal{C}$ and $\mathcal{D}$.
\begin{center}
    \begin{tikzpicture}

\draw  (-3,2) ellipse (0.6 and 1);
\draw  (-0.2,2) ellipse (0.6 and 1);
\node (v1) at (-3.2,1.6) {$c^\prime$};
\node (v2) at (-2.8,2.4) {$c$};
\node (v3) at (-0.4,2.4) {$d$};
\node (v4) at (0,1.6) {$d^\prime$};
\draw[->]  (v1) edge node[left, font = \scriptsize]{$f$} (v2);
\draw[->]  (v2) edge node[above, font = \scriptsize]{$x$} (v3);
\draw[->]  (v3) edge node[right, font = \scriptsize]{$g$} (v4);
\draw[->, dashed]  (v1) edge node[below, font = \scriptsize]{$g \circ x \circ f$} (v4);
\node at (-3.8,2.8) {$\mathcal{C}$};
\node at (0.6,2.8) {$\mathcal{D}$};
\end{tikzpicture}
\end{center}
Composition is prescribed by the functoriality of $u$. Given morphisms $f: c^\prime \rightarrow c$ in $\mathcal{C}$ and $g: d \rightarrow d^\prime$ in $\mathcal{D}$, so that $(f,g): (c, d) \rightarrow (c^\prime, d^\prime)$ is a morphism in $\mathcal{C}^{op} \times \mathcal{D}$, we define for $x \in u(c,d)$
$$g \circ x \circ f = u(f,g)(x) $$

There are interesting examples of collages. For a functor $F: \mathcal{C} \rightarrow \mathcal{D}$ the collage of the induced profunctor $F^*$ is precisely the mapping cylinder of $F$. $(\mathcal{C}, \mathcal{D})$-profunctors form a category with natural transformations as morphisms. The terminal object in this category assigns the singleton set to each pair of objects. Its collage gives us the join $\mathcal{C} * \mathcal{D}$. 

To make things even more interesting, observe that a collage $\textbf{col}(u)$ is naturally equipped with a map 
$$p: \textbf{col}(u) \rightarrow \Delta^1$$
where $\Delta^1 = \{0 \rightarrow 1 \}$ is the usual categorical 1-simplex, with $p^{-1}(0) \cong \mathcal{C}$ and $p^{-1}(0) \cong \mathcal{D}$. The following proposition due to B\'{e}nabou (\cite{street2001powerful}) tells us that any category equipped with a map to $\Delta^1$ is the collage of a profunctor. 

\begin{proposition}
\label{collage}
The collage construction produces an equivalence of categories between

\begin{itemize}
    \item $\textbf{Prof}_1$, the category whose objects are triples $(\mathcal{C}, \mathcal{D}, u)$ where $\mathcal{C}$, $\mathcal{D}$ are categories and $u$ is a profunctor between them, and a morphism $(F,G, \alpha) : (\mathcal{C}, \mathcal{D}, u) \rightarrow (\mathcal{C}^\prime, \mathcal{D}^\prime, u^\prime)$ between such triples consists of functors $F: \mathcal{C} \rightarrow \mathcal{C}^\prime$, $G: \mathcal{D} \rightarrow \mathcal{D}^\prime$, $\alpha: \textbf{col}(u) \rightarrow \textbf{col}(u^\prime)$ such that $\alpha|_\mathcal{C} = F$ and $\alpha|_\mathcal{D} = G$.
    
    \item The slice category $\textbf{Cat}/\Delta^1$, whose objects are pairs $(\mathcal{U},p)$ where $\mathcal{U}$ is a category and $p: \mathcal{U} \rightarrow \Delta^1$ is a functor, and whose morphisms between two such pairs $(\mathcal{U}, p)$ and $(\mathcal{V}, q)$ are functors $F: \mathcal{U} \rightarrow \mathcal{V}$ such that $q \circ F = p$. 
\end{itemize}
\end{proposition}

\begin{redi_proof}
It is clear that the collage construction extends into a functor 
$$\textbf{col} : \textbf{Prof}_1 \rightarrow \textbf{Cat}/\Delta^1$$
We may construct a weak inverse of $\textbf{col}$ in the obvious way. Let $p: \mathcal{U} \rightarrow \Delta^1$ be an object in $\textbf{Cat}/\Delta^1$ and let $\mathcal{U}_0 = p^{-1}(0)$ and $\mathcal{U}_1 = p^{-1}(1)$ be the corresponding fibers. Then we may define a $(\mathcal{U}_0, \mathcal{U}_1)$-profunctor $u_p$ simply by
$$u_p(x,y) = \mathcal{U}(x,y)$$
for $x \in \mathcal{U}_0$ and $y \in \mathcal{U}_1$. The functoriality of $u_p$ follows immediately. It is an easy exercise to verify that this construction is functorial and defines a weak inverse of $\textbf{col}$.
\begin{flushright}
$\square$
\end{flushright}
\end{redi_proof}

We would like to reserve the term \textit{correspondence} to indicate that we have maps to $\Delta^1$ in mind. Of course, all of the terms mentioned above are synonymous, at least for categories, but what they invoke feels (at least to us) different. 

As simple as Proposition \ref{collage} is, it allows us to transfer the concept of bimodule from the world of algebra into that of topology, because maps to $\Delta^1$ make sense in the category of simplicial sets and topological spaces. Being able to transfer ideas and concepts from one field to another is a known consequence of using the language of category theory. Moreover, this allows us to state and prove new results as well. For example in this work we will prove a classification theorem for inner fibrations (and simplicial maps in general). 

\subsection{Composition a.k.a. tensoring} \label{tensor}

We said profunctors are another notion of morphism between categories, therefore we should prescribe a composition. Given that profunctors are bimodules, their composition should be a tensor product. Let $u$ be a $(\mathcal{C}, \mathcal{D})$-profunctor and $v$ be a $(\mathcal{D}, \mathcal{E})$-profunctor. We define their composition to be ``the" $(\mathcal{C}, \mathcal{E})$-profunctor $v \otimes_\mathcal{D} u$ whose evaluation at a pair of objects $c \in \mathcal{C}$, $e \in \mathcal{E}$ is given by the coend formula
$$(v \otimes_\mathcal{D} u)(c,e) = \int_{d \in \mathcal{D}} v(d,e) \times u(c,d)$$

While this coend formula has the virtue of applying in any enrichment, it might not be very illuminating. It is useful to interpret the tensor product in terms of the collage construction. Consider $\textbf{col}(u)$ and $\textbf{col}(v)$ and juxtapose them along $\mathcal{D}$. What results is not a category because we cannot compose $x \in u(c,d)$ with $y \in v(d,e)$. Resolve this issue by generating a free category out of the data by declaring a new morphism $y \otimes x: c \rightarrow e$ serving as a composite . 
\begin{center}
	    \begin{tikzpicture}

\draw  (-3,2) ellipse (0.6 and 1);
\draw  (-0.4,3) ellipse (0.6 and 1);
\draw  (2.2,2) ellipse (0.6 and 1);

\node (v1) at (-3,2) {$c$};
\node (v2) at (-0.4,3) {$d$};
\node (v3) at (2.2,2) {$e$};
\draw[->, bend left = 5]  (v1) edge node[above, font = \scriptsize]{$x$} (v2);
\draw[->, bend left = 5]  (v2) edge node[above, font = \scriptsize]{$y$} (v3);
\draw[->, dashed, bend right = 10]  (v1) edge node[below, font = \scriptsize]{$y \otimes x$} (v3);

\node at (-3.9,2.4) {$\mathcal{C}$};
\node at (-1,4) {$\mathcal{D}$};
\node at (3.1,2.4) {$\mathcal{E}$};
\end{tikzpicture}
	\end{center}
This way a morphism in this new category between $c \in \mathcal{C}$ and $e \in \mathcal{E}$ is represented by a triple $(x,d,y)$ where $d \in \mathcal{D}$, $x \in u(c,d)$ and $y \in v(d,e)$. 
Then we impose the usual relation of the tensor product: for a morphism $f: d \rightarrow d^\prime$ in $\mathcal{D}$ we identify 
$$(fx, d^\prime, y) = (x, d, yf)$$
for all $x \in u(c,d)$ and $y \in v(d^\prime, e)$.
Finally we remove the objects of $\mathcal{D}$ to obtain the collage of $v \otimes_\mathcal{D} u$.
	
Equivalently, let $p: \mathcal{U} \rightarrow \Delta^1$ and $q: \mathcal{V} \rightarrow \Delta^1$ be composable correspondences, i.e. $p^{-1}(1) \cong \mathcal{D} \cong
q^{-1}(0)$ for some category $\mathcal{D}$. By taking the pushout along $\mathcal{D}$ we obtain a map $p \coprod_\mathcal{D} q : \mathcal{U} \coprod_\mathcal{D} \mathcal{V} \rightarrow \Delta^2 \cong \Delta^1 \coprod_{\Delta^0} \Delta^1$ (this is precisely the above picture). Then we produce the tensor product by taking the pullback
\begin{center}
    \begin{tikzpicture}

\node (v1) at (-3.5,3) {$\textbf{col}(v \otimes_\mathcal{D} u)$};
\node (v4) at (-3.5,1) {$\Delta^1$};
\node (v2) at (-0.5,3) {$\mathcal{U} \coprod_\mathcal{D} \mathcal{V}$};
\node (v3) at (-0.5,1) {$\Delta^2$};
\draw[->, dashed]  (v1) edge (v2);
\draw[->]  (v2) edge  (v3);
\draw[->, dashed]  (v1) edge  (v4);
\draw[->]  (v4) edge node[above, font = \scriptsize]{$d_1$} (v3);
\end{tikzpicture}
\end{center}
 
However, this composition operation does not produce a category because it is unital and associative only up to canonical isomorphism (we defined it using universal properties after all). Instead, we obtain a weak 2-category $\textbf{Cat}_{prof}$ whose objects are categories, 1-morphisms are profunctors and 2-morphisms are natural transformations between them. 

The unit profunctor for a category $\mathcal{C}$ is simply the profunctor $1_\mathcal{C}^*$ induced by the identity functor $1_\mathcal{C}$,
which assigns to a pair of objects $(c, c^\prime)$ the set of morphisms $\mathcal{C}(c, c^\prime) $. Its collage is the cylinder $\mathcal{C} \times \Delta^1$, seen as a correspondence via the projection map to $\Delta^1$.  

Now we observe the following. Consider a functor $F: \mathcal{X} \rightarrow \mathcal{C}$ between categories and let $f: a \rightarrow b$ be a morphism in $\mathcal{C}$. Let the categories $\mathcal{X}_a = F^{-1}(a)$ and $\mathcal{X}_b = F^{-1}(b)$ be the fibers of $a$ and $b$, and $\mathcal{X}_f = F^{-1}(f)$ the fibre of $f$, which is obtained by the pulling back $F$ along the functor $\Delta^1 \rightarrow \mathcal{C}$ which picks $f$,
\begin{center}
    \begin{tikzpicture}

\node (v1) at (-2.5,3) {$\mathcal{X}_f$};
\node (v4) at (-2.5,1) {$\Delta^1$};
\node (v2) at (-0.5,3) {$\mathcal{X}$};
\node (v3) at (-0.5,1) {$\mathcal{C}$};
\draw[->, dashed]  (v1) edge (v2);
\draw[->]  (v2) edge node[right, font = \scriptsize]{$F$} (v3);
\draw[->, dashed]  (v1) edge  (v4);
\draw[->]  (v4) edge node[above, font = \scriptsize]{$f$} (v3);
\end{tikzpicture}
\end{center}
is by definition an $(\mathcal{X}_a, \mathcal{X}_b)$-profunctor. 

Our definition of composition is fine enough to make the assignment $a \mapsto \mathcal{X}_a$ of fibers to each object in $\mathcal{C}$ functorial, so that each functor $F$ into $\mathcal{C}$ straightens into a diagram
$$\chi_F : \mathcal{C} \rightarrow \textbf{Cat}_{prof}$$
These diagrams turn out to be lax, meaning that if $f$ and $g$ are composable arrows in $\mathcal{C}$, we have comparison maps $\mathcal{X}_g \otimes \mathcal{X}_f \Rightarrow \mathcal{X}_{gf}$ which satisfy some coherence laws. 

In turn, we expect the assignment $F \mapsto \chi_F$ to define a functor
$$\chi : \textbf{Cat}/\mathcal{C} \rightarrow [\mathcal{C}, \textbf{Cat}_{prof}]_{lax}$$
from the slice category of functors over $\mathcal{C}$ into a category of lax $\mathcal{C}$-indexed diagrams of profunctors. Usual 2-categorical transformations between such functors are not the correct notion of morphism for our purposes though. To meet our expectations we need the notion of vertical transformation coming from double category theory.

It turns out $\chi$ is an equivalence of categories, so in some sense we have a classification result for functors. Quite interestingly, the weak inverse of $\chi$ is the double colimit functor defined and studied in \cite{grandis1999limits}. We will present these results in detail in Section \ref{pattern}.

\subsection{Double categories} \label{dblcat}

\begin{definition}
A (strict) double category is a category object in \textbf{Cat}.
\end{definition}
Let us unpack and depict the above definition. First, a category object in \textbf{Cat}, say $\mathbb{D}$, is comprised of a category $\mathbb{D}_0$ (the object category), a category $\mathbb{D}_1$ (the morphisms category) and source and target functors
$$s,t : \mathbb{D}_1 \rightarrow \mathbb{D}_0$$

We refer to the objects of $\mathbb{D}_0$ as the objects of $\mathbb{D}$. Each object in $u \in \mathbb{D}_1$ has a source and a target in $\mathbb{D}_0$ and hence it makes sense to depict them as arrows 
$$s(u) \xrightarrow{u} t(u)$$
Therefore we have two types of arrows between the objects of $\mathbb{D}$: those in $\mathbb{D}_0$ and the objects of $\mathbb{D}_1$. 
We will draw the arrows in $\mathbb{D}_0$ in the vertical direction and the objects of $\mathbb{D}_1$ as arrows in the horizontal direction. 

Since $s$ and $t$ are functors, a morphism $\alpha : u \rightarrow v$ in $\mathbb{D}_1$ has a source and target as well. We will depict these morphisms as square-shaped 2-cells
\begin{center}
        \begin{tikzpicture}

\node (v1) at (-2.5,3) {$s(u)$};
\node (v4) at (-2.5,1) {$s(v)$};
\node (v2) at (-0.5,3) {$t(u)$};
\node (v3) at (-0.5,1) {$t(v)$};
\draw[->]  (v1) edge node[above] {$u$} (v2);
\draw[->]  (v2) edge node[right] {$t(\alpha)$} (v3);
\draw[->]  (v1) edge node[left] {$s(\alpha)$} (v4);
\draw[->]  (v4) edge node[below] {$v$} (v3);
\node (v5) at (-1.5,2.6) {};
\node (v6) at (-1.5,1.4) {};
\draw[->]  (v5) edge[double] node[right] {$\alpha$} (v6);
\end{tikzpicture}
    \end{center}
We will refer to composition in $\mathbb{D}_1$ as vertical composition. 
\begin{center}
    \begin{tikzpicture}
\node (v1) at (-2.5,3) {$a_0$};
\node (v4) at (-2.5,1) {$b_0$};
\node (v2) at (-0.5,3) {$a_1$};
\node (v3) at (-0.5,1) {$b_1$};
\draw[->]  (v1) edge node[above, font = \scriptsize] {$u$} (v2);
\draw[->]  (v2) edge node[right, font = \scriptsize] {$f_1$} (v3);
\draw[->]  (v1) edge node[left, font = \scriptsize] {$f_0$} (v4);
\draw[->]  (v4) edge node[below, font = \scriptsize] {$v$} (v3);
\node (v5) at (-1.5,2.6) {};
\node (v6) at (-1.5,1.4) {};
\draw[->]  (v5) edge[double] node[right] {$\alpha$} (v6);

\node (v7) at (-2.5,-1) {$c_0$};
\node (v8) at (-0.5,-1) {$c_1$};
\draw[->]  (v4) edge node[left, font = \scriptsize] {$g_0$} (v7);
\draw[->]  (v3) edge node[right, font = \scriptsize] {$g_1$} (v8);
\draw[->]  (v7) edge node[below, font = \scriptsize] {$w$} (v8);
\node (v9) at (-1.5,0.5) {};
\node (v10) at (-1.5,-0.5) {};
\draw[->]  (v9) edge[double] node[right] {$\beta$} (v10);

\node (v11) at (0.5,1) {};
\node (v12) at (2,1) {};
\draw[|->]  (v11) edge node[above, font = \scriptsize] {compose} (v12);

\node (v13) at (3,2) {$a_0$};
\node (v14) at (5,2) {$a_1$};
\node (v15) at (3,0) {$c_0$};
\node (v16) at (5,0) {$c_1$};
\draw[->]  (v13) edge node[above, font = \scriptsize] {$u$} (v14);
\draw[->]  (v13) edge node[left, font = \scriptsize] {$g_0f_0$} (v15);
\draw[->]  (v14) edge node[right, font = \scriptsize] {$g_1f_1$} (v16);
\draw[->]  (v15) edge node[below, font = \scriptsize] {$w$} (v16);
\node (v17) at (4,1.5) {};
\node (v18) at (4,0.5) {};
\draw[->]  (v17) edge[double] node[right] {$\beta\alpha$} (v18);
\end{tikzpicture}
\end{center}

Another piece of data in a double category $\mathbb{D}$ is a composition functor 
$$\otimes : \mathbb{D}_1 \times_{\mathbb{D}_0} \mathbb{D}_1 \rightarrow \mathbb{D}_1$$
On objects, we interpret $\otimes$ as composition of horizontal morphisms in $\mathbb{D}$. On morphisms we interpret it as horizontal composition of 2-cells.

\begin{center}
    \begin{tikzpicture}
\node (v1) at (-2.5,3) {$a_0$};
\node (v4) at (-2.5,1) {$b_0$};
\node (v2) at (-0.5,3) {$a_1$};
\node (v3) at (-0.5,1) {$b_1$};
\draw[->]  (v1) edge node[above, font = \scriptsize] {$u_0$} (v2);
\draw[->]  (v2) edge node[right, font = \scriptsize] {$g$} (v3);
\draw[->]  (v1) edge node[left, font = \scriptsize] {$f$} (v4);
\draw[->]  (v4) edge node[below, font = \scriptsize] {$v_0$} (v3);
\node (v5) at (-1.5,2.6) {};
\node (v6) at (-1.5,1.4) {};
\draw[->]  (v5) edge[double] node[right] {$\alpha$} (v6);

\node (v7) at (1.5,3) {$a_2$};
\node (v8) at (1.5,1) {$b_2$};
\draw[->]  (v2) edge node[above, font = \scriptsize] {$u_1$} (v7);
\draw[->]  (v3) edge node[below, font = \scriptsize] {$v_1$} (v8);
\draw[->]  (v7) edge node[right, font = \scriptsize] {$h$} (v8);
\node (v9) at (0.5,2.5) {};
\node (v10) at (0.5,1.5) {};
\draw[->]  (v9) edge[double] node[right] {$\beta$} (v10);

\node (v11) at (2.5,2) {};
\node (v12) at (4,2) {};
\draw[|->]  (v11) edge node[above, font = \scriptsize] {compose} (v12);

\node (v13) at (5,3) {$a_0$};
\node (v14) at (7,3) {$a_2$};
\node (v15) at (5,1) {$b_0$};
\node (v16) at (7,1) {$b_2$};
\draw[->]  (v13) edge node[above, font = \scriptsize] {$u_1 \otimes u_0$} (v14);
\draw[->]  (v13) edge node[left, font = \scriptsize] {$f$} (v15);
\draw[->]  (v14) edge node[right, font = \scriptsize] {$h$} (v16);
\draw[->]  (v15) edge node[below, font = \scriptsize] {$v_1 \otimes v_0$} (v16);
\node (v17) at (6,2.5) {};
\node (v18) at (6,1.5) {};
\draw[->]  (v17) edge[double] node[right, font = \scriptsize] {$\beta \otimes \alpha$} (v18);
\end{tikzpicture}
\end{center}
Functionality of $\otimes$ translates into an interchange law, i.e. given the configuration of 2-cells 
\begin{center}
    \begin{tikzpicture}
\node (v1) at (-2.5,3) {$\bullet$};
\node (v4) at (-2.5,1) {$\bullet$};
\node (v2) at (-0.5,3) {$\bullet$};
\node (v3) at (-0.5,1) {$\bullet$};
\draw[->]  (v1) edge  (v2);
\draw[->]  (v2) edge  (v3);
\draw[->]  (v1) edge  (v4);
\draw[->]  (v4) edge  (v3);
\node (v5) at (-1.5,2.6) {};
\node (v6) at (-1.5,1.4) {};
\draw[->]  (v5) edge[double] node[right] {$\alpha$} (v6);

\node (v7) at (1.5,3) {$\bullet$};
\node (v8) at (1.5,1) {$\bullet$};
\draw[->]  (v2) edge  (v7);
\draw[->]  (v3) edge  (v8);
\draw[->]  (v7) edge  (v8);
\node (v9) at (0.5,2.5) {};
\node (v10) at (0.5,1.5) {};
\draw[->]  (v9) edge[double] node[right] {$\beta$} (v10);

\node (v11) at (-2.5,-1) {$\bullet$};
\node (v12) at (-0.5,-1) {$\bullet$};
\node (v13) at (1.5,-1) {$\bullet$};
\draw[->]  (v4) edge (v11);
\draw[->]  (v3) edge (v12);
\draw[->]  (v8) edge (v13);
\draw[->]  (v11) edge (v12);
\draw[->]  (v12) edge (v13);
\node (v14) at (-1.5,0.5) {};
\node (v15) at (-1.5,-0.5) {};
\draw[->]  (v14) edge[double] node[right] {$\gamma$} (v15);
\node (v16) at (0.5,0.5) {};
\node (v17) at (0.5,-0.5) {};
\draw[->]  (v16) edge[double] node[right] {$\theta$} (v17);
\end{tikzpicture}
\end{center}
the two different ways of obtaining a total composite (first composing horizontally and then vertically and vice versa) coincide. Equationally this is written as
$$(\beta \otimes \alpha)(\theta \otimes \gamma) = (\theta\beta) \otimes (\gamma\alpha)$$

$\otimes$ is associative and unital. We will depict unit horizontal morphisms as double edges
\begin{center}
    \begin{tikzpicture}

\node (v1) at (-3,1) {$\bullet$};
\node (v2) at (-1,1) {$\bullet$};
\draw[-]  (v1) edge[double] (v2);
\end{tikzpicture}
\end{center}
and the unit vertical and horizontal 2-cells accordingly. 

There are interesting examples of strict double categories but since they are two-dimensional categorical structures weak versions are more useful. For our purposes we are interested in double categories which are weak (unital and associative only up to coherent isomorphism) in the horizontal direction. This means we postulate some extra pieces of data:
\begin{itemize}
    \item an invertible associator 2-cell for each composable triple $a \xrightarrow{u} b \xrightarrow{v} c \xrightarrow{w} d$ of horizontal morphisms
\begin{center}
    \begin{tikzpicture}
\node (v1) at (-2.5,3) {$a$};
\node (v4) at (-2.5,1) {$a$};
\node (v2) at (-0.5,3) {$c$};
\node (v3) at (-0.5,1) {$b$};
\draw[->]  (v1) edge node[above, font = \scriptsize] {$v \otimes u$} (v2);

\draw[->]  (v1) edge node[left, font = \scriptsize] {$1_a$} (v4);
\draw[->]  (v4) edge node[below, font = \scriptsize] {$u$} (v3);
\node (v5) at (1.5,3) {$d$};
\node (v6) at (1.5,1) {$d$};
\draw[->]  (v2) edge node[above, font = \scriptsize] {$w$}  (v5);
\draw[->]  (v3) edge node[below, font = \scriptsize] {$w \otimes v$}  (v6);
\draw[->]  (v5) edge node[right, font = \scriptsize] {$1_d$} (v6);
\node (v7) at (-0.5,2.5) {};
\node (v8) at (-0.5,1.5) {};
\draw[->]  (v7) edge[double] node[right, font = \scriptsize] {$\cong$} (v8);
\end{tikzpicture}
\end{center}

\item invertible unitor 2-cells for each horizontal morphism $u: a \rightarrow b$
\begin{center}
    \begin{tikzpicture}

\node (v1) at (-3,1) {$a$};
\node (v2) at (-1.5,1) {$a$};
\node (v3) at (0,1) {$b$};
\node (v4) at (-3,-1) {$a$};
\node (v5) at (0,-1) {$b$};
\draw[-]  (v1) edge[double] (v2);
\draw[->]  (v2) edge node[above, font = \scriptsize] {$u$} (v3);
\draw[->]  (v1) edge node[left, font = \scriptsize] {$1_a$} (v4);
\draw[->]  (v3) edge node[right, font = \scriptsize] {$1_b$} (v5);
\draw[->]  (v4) edge node[below, font = \scriptsize] {$u$} (v5);
\node (v6) at (-1.5,0.5) {};
\node (v7) at (-1.5,-0.5) {};
\draw[->]  (v6) edge[double] node[right, font = \scriptsize] {$\cong$} (v7);
\end{tikzpicture} \hspace{1 cm}
 \begin{tikzpicture}

\node (v1) at (-3,1) {$a$};
\node (v2) at (-1.5,1) {$b$};
\node (v3) at (0,1) {$b$};
\node (v4) at (-3,-1) {$a$};
\node (v5) at (0,-1) {$b$};
\draw[->]  (v1) edge node[above, font = \scriptsize] {$w$} (v2);
\draw[-]  (v2) edge[double] (v3);
\draw[->]  (v1) edge node[left, font = \scriptsize] {$1_a$} (v4);
\draw[->]  (v3) edge node[right, font = \scriptsize] {$1_b$} (v5);
\draw[->]  (v4) edge node[below, font = \scriptsize] {$u$} (v5);
\node (v6) at (-1.5,0.5) {};
\node (v7) at (-1.5,-0.5) {};
\draw[->]  (v6) edge[double] node[right, font = \scriptsize] {$\cong$} (v7);
\end{tikzpicture}
\end{center}
\end{itemize}
satisfying the usual 2-categorical coherence laws (see Definition 5.2.1 in \cite{leinster2004higher}). 

\begin{example}
There is a (weak) double category \textbf{Prof} whose
\begin{itemize}
    \item objects are categories
    \item vertical morphisms are functors
    \item horizontal morphisms are profunctors
    \item 2-cells of the form 
    \begin{center}
        \begin{tikzpicture}
\node (v1) at (-2.5,3) {$\mathcal{C}_0$};
\node (v4) at (-2.5,1) {$\mathcal{D}_0$};
\node (v2) at (-0.5,3) {$\mathcal{C}_1$};
\node (v3) at (-0.5,1) {$\mathcal{D}_1$};
\draw[->]  (v1) edge node[above, font = \scriptsize] {$u$} (v2);
\draw[->]  (v2) edge node[right, font = \scriptsize] {$G$} (v3);
\draw[->]  (v1) edge node[left, font = \scriptsize] {$F$} (v4);
\draw[->]  (v4) edge node[below, font = \scriptsize] {$v$}  (v3);
\node (v5) at (-1.5,2.6) {};
\node (v6) at (-1.5,1.4) {};
\draw[->]  (v5) edge[double] node[right] {$\alpha$} (v6);
        \end{tikzpicture}
    \end{center}
are functors $\alpha: \textbf{col}(u) \rightarrow \textbf{col}(v)$ such that $\alpha|_{\mathcal{C}_0} = F$ and $\alpha|_{\mathcal{C}_1} = G$

\item vertical composition is just composition of functors and horizontal composition is the tensor product
\end{itemize}
\end{example}

In a similar fashion one can construct all sorts of double categories of bimodules and more (see \cite{leinster2004higher} or \cite{grandis1999limits}). Here \textbf{Prof} is our main object of interest, as it incorporates in it both functors and profunctors. 

A double category $\mathbb{D}$ has incorporated in it two 2-categories, one horizontal and one vertical. The vertical 2-category $\mathbb{D}_v$ has as 2-morphisms 2-cells in $\mathbb{D}$ of the form
\begin{center}
        \begin{tikzpicture}

\node (v1) at (-2.5,3) {$a$};
\node (v4) at (-2.5,1) {$b$};
\node (v2) at (-0.5,3) {$a$};
\node (v3) at (-0.5,1) {$b$};
\draw[-]  (v1) edge[double]  (v2);
\draw[->]  (v2) edge node[right] {$g$} (v3);
\draw[->]  (v1) edge node[left] {$f$} (v4);
\draw[-]  (v4) edge[double] (v3);
\node (v5) at (-1.5,2.6) {};
\node (v6) at (-1.5,1.4) {};
\draw[->]  (v5) edge[double] node[right] {$\alpha$} (v6);
\end{tikzpicture}
    \end{center}
and the horizontal 2-category $\mathbb{D}_h$ has as 2-morphisms 2-cells in $\mathbb{D}$ of the form
\begin{center}
        \begin{tikzpicture}

\node (v1) at (-2.5,3) {$a$};
\node (v4) at (-2.5,1) {$a$};
\node (v2) at (-0.5,3) {$b$};
\node (v3) at (-0.5,1) {$a$};
\draw[->]  (v1) edge node[above] {$u$} (v2);
\draw[->]  (v2) edge node[right] {$1_b$} (v3);
\draw[->]  (v1) edge node[left] {$1_a$} (v4);
\draw[->]  (v4) edge node[below] {$v$} (v3);
\node (v5) at (-1.5,2.6) {};
\node (v6) at (-1.5,1.4) {};
\draw[->]  (v5) edge[double] node[right] {$\alpha$} (v6);
\end{tikzpicture}
    \end{center}
If $\mathbb{D}$ is horizontally weak then $\mathbb{D}_h$ is weak. 
$\textbf{Prof}_v$ is the usual 2-category of categories, functors and natural transformations, while $\textbf{Prof}_h$ is simply $\textbf{Cat}_{prof}$. 

Further double categorical concepts (horizontal diagrams, vertical transformations and double colimits) will be presented as needed in Section \ref{pattern}.
We will have a taste throughout this work of how powerful of an organizing principle double categories are and how efficient their language is. 

\subsection{The meaning of inner fibrations} \label{inner}

Let $\wedge^n_i$ be the simplicial subset of $\Delta^n$ generated by all the faces except the $i$-th one, for $[n] \in \Delta$. Recall that a map of simplicial sets $p: X \rightarrow A$ is called an inner fibration if it satisfies the left lifting property with respect to all inclusions $\wedge^n_i \subset \Delta^n$ for all $n$ and $0 < i < n$, i.e. given the solid square 
\begin{center}
    \begin{tikzpicture}

\node (v1) at (-4,3.5) {$\wedge^n_i$};
\node (v3) at (-4,1.5) {$\Delta^n$};
\node (v2) at (-2,3.5) {$X$};
\node (v4) at (-2,1.5) {$A$};
\draw[->]  (v1) edge (v2);
\draw[left hook ->]  (v1) edge (v3);
\draw[->]  (v2) edge node[right, font = \scriptsize]{$p$} (v4);
\draw[->]  (v3) edge (v4);
\draw[->, dashed]  (v3) edge (v2);
\end{tikzpicture}
\end{center}
there is a dotted arrow making the diagram commute. 

Inner fibrations are important objects of study in the theory of higher categories. For example, all fibrations in the Joyal model structure (\cite[Section 2.2.5]{lurie2009higher}) on simplicial sets are inner fibrations. Also recall that a simplicial set $A$ is called an $\infty$-category if the terminal map $A \rightarrow \Delta^0$ is an inner fibration. 

Let $p: X \rightarrow A$ be an inner fibration between simplicial sets. For an $n$-simplex $\sigma \in A_n$ consider the fiber $X_\sigma = p^{-1}(\sigma)$ obtained by the pullback square
\begin{center}
    \begin{tikzpicture}

\node (v1) at (-2.5,3) {$X_\sigma$};
\node (v4) at (-2.5,1) {$\Delta^n$};
\node (v2) at (-0.5,3) {$X$};
\node (v3) at (-0.5,1) {$A$};
\draw[->, dashed]  (v1) edge (v2);
\draw[->]  (v2) edge node[right, font = \scriptsize]{$p$} (v3);
\draw[->, dashed]  (v1) edge node[left, font = \scriptsize]{$p_\sigma$} (v4);
\draw[->]  (v4) edge node[above, font = \scriptsize]{$\sigma$} (v3);
\end{tikzpicture}
\end{center}
Fibrations are stable under pullback, so $p_\sigma$ is an inner fibration as well. 

Moreover, if the target of an inner fibration is a category then the source is an $\infty$-category. In our case we conclude that fibers $p^{-1}(\sigma)$ over each simplex $\sigma$ have to be $\infty$-categories. It is not difficult to see that the converse is also true: if a map $p: X \rightarrow A$ of simplicial sets is such that fibers over each simplex are $\infty$-categories then it is an inner fibration. 

Referring to the above observations we quote the authors of \cite{barwick2018fibrations}
\begin{quote}
    ``So in a strong sense, we'll understand the ``meaning" of inner fibrations once we understand the ``meaning" of functors from $\infty$-categories to $\Delta^n$"
\end{quote}
In analogy with category theory maps from an $\infty$-category to $\Delta^1$ can be called correspondences. In this vein, maps to $\Delta^n$ deserve to be called \textit{higher correspondences}.

A higher correspondence $p : X \rightarrow \Delta^n$ may be thought of as consisting of $(n+1)$ $\infty$-categories $X_i = p^{-1}(i)$, $i = 0, 1 \dots n$, and a collage $X$ between them. By the latter we mean that $X$ is formed from the $X_i$'s by ``adding" 1-simplices which join vertices of $X_i$ and $X_j$ only if $i<j$, and higher simplices after that. 

Given an inner fibration $p : X \rightarrow A $, for each simplex $\sigma \in A_n$ we obtain a correspondence $X_\sigma$. We would like to see the assignment 
$$\sigma \mapsto X_\sigma$$
as being ``functorial". As pointed out by Lurie (\cite[s.2.3.1]{lurie2009higher}), in order for functoriality to make sense here there should exist some higher category, say \textbf{Corr}, whose objects are $\infty$-categories, 1-morphisms are correspondences and higher morphisms are higher correspondences. 

So, once \textbf{Corr} is defined correctly each inner fibration $p$ with target a simplicial set $A$ should produce a diagram 
$$\chi_p : A \rightarrow \textbf{Corr}$$
Moreover $p \mapsto \chi_p$ is also expected to be functorial and to define some sort of equivalence
$$\chi: (\textbf{sSet}/A)_{inner} \xrightarrow{\sim} [A, \textbf{Corr}]$$
between inner fibrations over $A$ and $A$-indexed diagrams of correspondences. 

A natural approach in achieving \textbf{Corr} would be to understand correspondences as $\infty$-categorical analogs of profunctors. This would consist of replacing the category \textbf{Set} by the $\infty$-category of spaces, and the coend formula by a homotopical analog. Such profunctors are equivalent to maps from $\infty$-categories to $\Delta^1$ (\cite{stevenson2018model}). 

For ordinary categories a functor to $\Delta^2$ is the same as the data for profunctors $u$, $v$, $w$ and a natural transformation $v \otimes u \Rightarrow w$ (we come back to this in detail in the next section). A similar statement is expected to hold for $\infty$-categories. This indicates that if \textbf{Corr} is to be realized as an $\infty$-category, then extra conditions are required to ensure that transformations such as above are invertible. Such an $\infty$-category is constructed in \cite{ayala2020fibrations}, and it serves as a classifier for certain fibrations called exponentiable, but not inner fibrations.

Inner fibrations aside, the above discussion extends to all simplicial maps. If $f: X \rightarrow A$ is a map of simplicial sets then it is expected of the assignment $\sigma \mapsto X_\sigma$, for simplices $\sigma \in A$, to provide a straightening of $f$ into a map from $A$ into a higher structure of simplicial sets and correspondences of different dimensions between them. 

Our strategy is to define a version of \textbf{Corr} for simplicial sets and prove that all simplicial maps straighten as expected. Then we will obtain a classifying result for inner fibrations as a corollary. We will obtain $\chi$ as an equivalence of categories and sketch how to promote it to an equivalence of \textbf{sSet}-categories in \ref{class_results}.

\section[The fibrational perspective]{A pattern: functions and functors from the fibrational  perspective} \label{pattern}

We would like to study the relationship between higher correspondences and simplicial maps (or inner fibrations) as part of a larger pattern. This pattern may be seen, on a fundamental level, to begin with sets and functions. 

Let $f: X \rightarrow A$ be a function between two sets. We observe that all the information about $f$ lies in its fibers in the following way. Let $X_a = f^{-1}(a)$, $a \in A$, be the fibers of $f$. The assignment $a \mapsto X_a$ defines a functor $$\chi_f : A \rightarrow \textbf{Set}$$
and $f \mapsto \chi_f$ itself produces a functor
$$\chi : \textbf{Set}/A \rightarrow [A, \textbf{Set}]$$
from the slice category of maps over $A$ into the category of $A$-indexed families of sets with natural transformations as morphisms.

Conversely, given $\phi \in [A, \textbf{Set}]$ we may form the disjoint union $\coprod_{a \in A} \phi(a)$ which is naturally equipped with a map to $A$ so that we obtain a functor
$$\coprod : [A, \textbf{Set}] \rightarrow \textbf{Set}/A$$
It is easy to see that $\chi$ and $\coprod$ are weak inverses of each other and they form an equivalence of categories. For instance, given $(X,f) \in \textbf{Set}/A$ we have a natural isomorphism $X \cong \coprod_{a \in A} X_a$. 

The discussion immediately becomes more complex when we study categories and functors in this fashion. Let $F: \mathcal{X} \rightarrow \mathcal{A}$ be a functor and let $\mathcal{X}_a = F^{-1}(a)$, $a \in \mathcal{A}$, be the fibers of $F$, i.e. $\mathcal{X}_a$ has as objects those objects of $\mathcal{X}$ which map to $a \in \mathcal{A}$ and as morphisms those morphisms in $\mathcal{X}$ which map to the identity $1_a$. $a \mapsto \mathcal{X}_a$ gives us a family of categories indexed by the objects of $\mathcal{A}$, but unlike sets, categories have morphisms, so this much information is far from being enough to reconstruct $(\mathcal{X}, F)$ up to isomorphism. 

Therefore we have to consider not only fibers over objects but fibers over morphisms as well if we hope to classify functors fibrationally. Let $f: a \rightarrow b$ be a morphism in $\mathcal{A}$. Let $\mathcal{X}_f = F^{-1}(f)$ be defined by pulling back along the map $\Delta^1 \rightarrow \mathcal{A}$ which picks $f$
\begin{center}
    \begin{tikzpicture}

\node (v1) at (-2.5,3) {$\mathcal{X}_f$};
\node (v4) at (-2.5,1) {$\Delta^1$};
\node (v2) at (-0.5,3) {$\mathcal{X}$};
\node (v3) at (-0.5,1) {$\mathcal{A}$};
\draw[->, dashed]  (v1) edge (v2);
\draw[->]  (v2) edge node[right, font = \scriptsize]{$F$} (v3);
\draw[->, dashed]  (v1) edge  (v4);
\draw[->]  (v4) edge node[above, font = \scriptsize]{$f$} (v3);
\end{tikzpicture}
\end{center}
Proposition \ref{collage} tells us that $\mathcal{X}_f$ is an $(\mathcal{X}_a, \mathcal{X}_b)$-profunctor. Let us see that this way $a \mapsto \mathcal{X}_a$ defines a lax 2-functor
$$\chi_F : \mathcal{A} \rightarrow \textbf{Cat}_{prof}$$
with $\chi_F(f) = \mathcal{X}_f$ for all morphisms $f$ in $\mathcal{A}$. 

Let $a \xrightarrow{f} b \xrightarrow{g} c$ be a pair of composable morphisms in $\mathcal{A}$. We may define the structure natural transformations $\phi_{f,g} : \mathcal{X}_g \otimes_{\mathcal{X}_b} \mathcal{X}_f \Rightarrow \mathcal{X}_{gf}$ in the obvious way as follows. Let $x \in \mathcal{X}_a$ and $y \in \mathcal{X}_c$ be two objects. Recall that an element of the set $(\mathcal{X}_g \otimes_{\mathcal{X}_b} \mathcal{X}_f) (x,y)$ is represented by a triple $(\alpha, z, \beta)$ where $z$ is an object of $\mathcal{X}_b$, $\alpha : x \rightarrow z$ is a morphism in $\mathcal{X}_f$ and  $\beta : z \rightarrow y$ is a morphism in $\mathcal{X}_g$. For such an element define 
$$\phi_{f,g}(\alpha, z, \beta) = \beta \circ \alpha$$
where the composition on the right is in $\mathcal{X}$. 

It is easy to see that this map is well-defined and that the above transformations make of $\chi_F$ a unitally strict lax 2-functor. By unitally strict we mean $\chi_F (1_a) = 1_{\mathcal{X}_a}^*$, which holds by construction. For this reason we assume all lax diagrams to be unitally strict throughout this section.

It is reasonable (and necessary to progress with results here) to expect the construction $F \mapsto \chi_F$ to be functorial, i.e. there should be a functor 
$$\chi : \textbf{Cat}/\mathcal{A} \rightarrow [\mathcal{A}, \textbf{Cat}_{prof}]_{lax}$$
from the slice category of functors over $\mathcal{A}$ to the category of lax $\mathcal{A}$-indexed diagrams of profunctors. However, usual (2-categorical) transformations between lax functors do not do the job for us (see also the Proposition in \cite{street2001powerful}). 

The correct notion of morphism in this case is provided by double category theory. Recall that $\textbf{Cat}_{prof}$ is the horizontal 2-category of the double category \textbf{Prof}, so we may interpret our functors $\chi_F$ for various $F$ to hit the horizontal part of \textbf{Prof}. The natural double categorical notion of morphism between such functors is that of a vertical transformation, which we present now. 

\begin{definition} \label{vtrans}
Let $\mathbb{D}$ be a double category, $\mathcal{A}$ a category and $(F, \phi), (G, \psi) : \mathcal{A} \rightarrow \mathbb{D}_h$ be two lax functors valued in the horizontal part of $\mathbb{D}$. Then a vertical natural transformation $\alpha : F \Rightarrow G$ consists of:
\begin{itemize}
    \item a vertical morphism $\alpha_a : F(a) \rightarrow G(a)$ for all $a \in \mathcal{A}$
    \item a 2-cell $\alpha_f$ as below for all morphisms $f: a \rightarrow b$ in $\mathcal{A}$
    \begin{center}
        \begin{tikzpicture}

\node (v1) at (-2.5,3) {$F(a)$};
\node (v4) at (-2.5,1) {$G(a)$};
\node (v2) at (-0.5,3) {$F(b)$};
\node (v3) at (-0.5,1) {$G(b)$};
\draw[->]  (v1) edge node[above] {$Ff$} (v2);
\draw[->]  (v2) edge node[right] {$\alpha_a$} (v3);
\draw[->]  (v1) edge node[left] {$\alpha_b$} (v4);
\draw[->]  (v4) edge node[below] {$Gf$} (v3);
\node (v5) at (-1.5,2.6) {};
\node (v6) at (-1.5,1.4) {};
\draw[->]  (v5) edge[double] node[right] {$\alpha_f$} (v6);
\end{tikzpicture}
    \end{center}
\end{itemize}
such that ``everything commutes", i.e. for all composable pairs of morphisms $a \xrightarrow{f} b \xrightarrow{g} c$ in $\mathcal{A}$ the following diagramatic equation holds
\begin{center}
    \begin{tikzpicture}

\node (v1) at (-4,3) {$Fa$};
\node (v2) at (-2,3) {$Fb$};
\node (v4) at (-4,1) {$Ga$};
\node (v5) at (-2,1) {$Gb$};
\node (v3) at (0,3) {$Fc$};
\node (v6) at (0,1) {$Gc$};
\node (v7) at (-4,-1) {$Gb$};
\node (v8) at (0,-1) {$Gc$};
\node at (1,1) {$=$};
\node (v9) at (2,3) {$Fa$};
\node (v10) at (4,3) {$Fb$};
\node (v11) at (6,3) {$Fc$};
\node (v12) at (2,1) {$Fa$};
\node (v13) at (6,1) {$Fc$};
\node (v14) at (2,-1) {$Ga$};
\node (v15) at (6,-1) {$Gc$};
\draw[->]  (v1) edge node[above, font = \scriptsize] {$Ff$} (v2);
\draw[->]  (v2) edge node[above, font = \scriptsize] {$Fg$} (v3);
\draw[->]  (v4) edge node[below, font = \scriptsize] {$Gf$} (v5);
\draw[->]  (v5) edge node[below, font = \scriptsize] {$Gg$} (v6);
\draw[->]  (v7) edge node[below, font = \scriptsize] {$G(gf)$} (v8);
\draw[->]  (v1) edge node[left, font = \scriptsize] {$\alpha_a$} (v4);
\draw[->]  (v2) edge node[right, font = \scriptsize] {$\alpha_b$} (v5);
\draw[->]  (v3) edge node[right, font = \scriptsize] {$\alpha_c$} (v6);
\draw[->]  (v4) edge node[left, font = \scriptsize] {$1_{Ga}$} (v7);
\draw[->]  (v6) edge  node[right, font = \scriptsize] {$1_{Gc}$} (v8);
\draw[->]  (v9) edge node[above, font = \scriptsize] {$Ff$} (v10);
\draw[->]  (v10) edge node[above, font = \scriptsize] {$Fg$} (v11);
\draw[->]  (v12) edge node[below, font = \scriptsize] {$F(gf)$} (v13);
\draw[->]  (v14) edge node[below, font = \scriptsize] {$G(gf)$} (v15);
\draw[->]  (v9) edge node[left, font = \scriptsize] {$1_{Fa}$} (v12);
\draw[->]  (v11) edge node[right, font = \scriptsize] {$1_{Fc}$} (v13);
\draw[->]  (v12) edge node[left, font = \scriptsize] {$\alpha_a$} (v14);
\draw[->]  (v13) edge node[right, font = \scriptsize] {$\alpha_c$} (v15);

\node (v16) at (-3,2.5) {};
\node (v17) at (-3,1.5) {};
\draw[->]  (v16) edge[double] node[right, font = \scriptsize] {$\alpha_f$} (v17);
\node (v18) at (-1,2.5) {};
\node (v19) at (-1,1.5) {};
\draw[->]  (v18) edge[double] node[right, font = \scriptsize] {$\alpha_g$} (v19);
\node (v20) at (-2,0.5) {};
\node (v21) at (-2,-0.5) {};
\draw[->]  (v20) edge[double] node[right, font = \scriptsize] {$\psi_{f,g}$} (v21);
\node (v22) at (4,2.5) {};
\node (v23) at (4,1.5) {};
\draw[->]  (v22) edge[double] node[right, font = \scriptsize] {$\phi_{f,g}$} (v23);
\node (v24) at (4,0.5) {};
\node (v25) at (4,-0.5) {};
\draw[->]  (v24) edge[double] node[right, font = \scriptsize] {$\alpha_{gf}$} (v25);

\end{tikzpicture}
\end{center}
or simply
$$\psi_{f,g} \circ (\alpha_g \otimes \alpha_f) = \alpha_{gf} \circ \phi_{f,g}$$
Here $\phi_{f,g}: Fg \otimes Ff \Rightarrow F(gf)$ and $\psi_{f,g} : Gg \otimes Gf \Rightarrow G(gf)$ in $\mathbb{D}_h$ are the structure maps corresponding to the composable pair $(f,g)$. 
\end{definition}

Such transformations compose by composing the structure vertical maps and cells, so that in general we have a category $[\mathcal{A}, \mathbb{D}_{h}]_{lax}$ of lax $\mathcal{A}$-indexed horizontal diagrams in $\mathbb{D}$. In particular, with this notion of morphism we obtain our functor
$$\chi : \textbf{Cat}/\mathcal{A} \rightarrow [\mathcal{A}, \textbf{Cat}_{prof}]_{lax}$$
Indeed, let $(\mathcal{X}, F), (\mathcal{Y}, G) \in \textbf{Cat}/\mathcal{A}$ and $U: \mathcal{X} \rightarrow \mathcal{Y}$ a morphism, i.e. $U$ is a functor such that the triangle 
\begin{center}
    \begin{tikzpicture}

\node (v1) at (-3,2) {$\mathcal{X}$};
\node (v2) at (-1,2) {$\mathcal{Y}$};
\node (v3) at (-2,0.6) {$\mathcal{A}$};
\draw[->]  (v1) edge node[above, font = \scriptsize] {$U$} (v2);
\draw[->]  (v1) edge node[left, font = \scriptsize] {$F$} (v3);
\draw[->]  (v2) edge node[right, font = \scriptsize] {$G$} (v3);
\end{tikzpicture}
\end{center}
commutes.
Then we have an induced vertical transformation $\chi(U) : \chi_F \Rightarrow \chi_G$ with 
\begin{itemize}
    \item structure vertical maps (functors) for $a \in \mathcal{A}$ the restrictions $U|_{\mathcal{X}_a}: \mathcal{X}_a \rightarrow \mathcal{Y}_a$ 
    \item structure 2-cells for morphisms $f: a \rightarrow b$ the restrictions $U|_{\mathcal{X}_f} : \mathcal{X}_f \rightarrow \mathcal{Y}_f$
    \begin{center}
        \begin{tikzpicture}

\node (v1) at (-2.5,3) {$\mathcal{X}_a$};
\node (v4) at (-2.5,1) {$\mathcal{Y}_a$};
\node (v2) at (-0.5,3) {$\mathcal{X}_b$};
\node (v3) at (-0.5,1) {$\mathcal{Y}_b$};
\draw[->]  (v1) edge node[above, font = \scriptsize] {$\mathcal{X}_f$} (v2);
\draw[->]  (v2) edge node[right, font = \scriptsize] {$U|_{\mathcal{X}_a}$} (v3);
\draw[->]  (v1) edge node[left, font = \scriptsize] {$U|_{\mathcal{X}_b}$} (v4);
\draw[->]  (v4) edge node[below, font = \scriptsize] {$\mathcal{Y}_f$} (v3);
\node (v5) at (-1.5,2.6) {};
\node (v6) at (-1.5,1.4) {};
\draw[->]  (v5) edge[double] node[right, font = \scriptsize] {$U|_{\mathcal{X}_f}$} (v6);
\end{tikzpicture}
    \end{center}
\end{itemize}

Lastly, we expect of $\chi$ to be an equivalence of categories, so that we interpret it as a classifying functor. In the case of sets, the inverse of $\chi$ is given by the colimiting process of the disjoint union. It turns out that the colimiting process involved in the case of functors is that of double colimit introduced and studied in \cite{grandis1999limits}. 

Again let $\mathbb{D}$ be a double category and 
$$F : \mathcal{A} \rightarrow \mathbb{D}_h$$
be a lax horizontal diagram indexed by a category $\mathcal{A}$. We define the category of vertical cones over $F$, $\textbf{vCone}(F)$, with:
\begin{itemize}
    \item objects pairs $(c, \alpha)$ where $c$ is an object of $\mathbb{D}$ and $\alpha : F \Rightarrow c$ is a vertical transformation from $F$ into the constant diagram at $c$ (which we denote simply $c$ by abusing notation)
    \item morphisms $ (c, \alpha) \rightarrow (d, \beta)$ between two such pairs being vertical morphisms $f: c \rightarrow d$ in $\mathcal{D}$ such that $\beta f = \alpha$
    \begin{center}
    \begin{tikzpicture}

\node (v1) at (-3,2) {$F$};
\node (v2) at (-1,2) {$c$};
\node (v3) at (-2,0.6) {$d$};
\draw[->]  (v1) edge[double] node[above, font = \scriptsize] {$\alpha$} (v2);
\draw[->]  (v1) edge[double] node[left, font = \scriptsize] {$\beta$} (v3);
\draw[->]  (v2) edge node[right, font = \scriptsize] {$f$} (v3);

\end{tikzpicture}
\end{center}
In the diagram above $f$ is to be understood as the vertical transformation induced by the vertical morphism $f$.
    
\end{itemize}
The special feature to be noted in this category is that we require morphisms between vertical cones to come from vertical morphisms and not simply be vertical transformations between constant functors. 

\begin{definition}
In a setting as above, the double colimit of the functor $F$ is an object $\textbf{dcolim}F \in \mathbb{D}$ which is initial in the category of vertical cones $\textbf{vCone}(F)$
\end{definition}
Clearly, double colimits are functorial and if they exist they are unique up to unique isomorphism in $\mathbb{D}_0$. The following proposition (see also Example 4.1 in \cite{niefield2012glueing}) provides a model for all double colimits in \textbf{Prof}.

\begin{proposition} \label{model}
Let $\mathcal{A}$ be a category and 
$$(F, \phi) : \mathcal{A} \rightarrow \textbf{Prof}_h$$
be a (strictly unital) lax functor. Let $\mathcal{X}$ be the category defined with:
\begin{itemize}
    \item objects pairs $(a,x)$ where $a \in \mathcal{A}$, $x \in F(a)$
    \item morphisms $(f, u) : (a,x) \rightarrow (b,y)$ between such pairs given by $f: a \rightarrow b$ morphism in $\mathcal{A}$ and $u \in (Ff)(x,y)$
    \item composition defined by the formula
    $$(g,v) \circ (f,u) = (gf, \phi(y \otimes x))$$
\end{itemize}
Then 
$$\mathcal{X} \cong \textbf{dcolim}(F)$$
\end{proposition}

\begin{redi_proof}
First, the coherence laws for lax functors and the relations in the definition of the tensor product make sure that $\mathcal{X}$ is well-defined. We will prove this category is a model for the double colimit of $F$ by direct verification of the universal property. 

Define the structure vertical transformation $\alpha : F \Rightarrow \mathcal{X}$ with
\begin{itemize}
    \item structure functors $\alpha_a : F(a) \rightarrow \mathcal{X}$, $a \in \mathcal{A}$, given by $\alpha_a(x) = (a,x)$ on objects and $\alpha_a(f) = (1_a, f)$ on morphisms
    \item structure 2-cells, for morphisms $f: a \rightarrow b$ in $\mathcal{A}$, 
    \begin{center}
        \begin{tikzpicture}

\node (v1) at (-2.5,3) {$F(a)$};
\node (v4) at (-2.5,1) {$\mathcal{X}$};
\node (v2) at (-0.5,3) {$F(b)$};
\node (v3) at (-0.5,1) {$\mathcal{X}$};
\draw[->]  (v1) edge node[above, font = \scriptsize] {$Ff$} (v2);
\draw[->]  (v2) edge node[right, font = \scriptsize] {$\alpha_b$}  (v3);
\draw[->]  (v1) edge node[left, font = \scriptsize] {$\alpha_a$} (v4);
\draw[-]  (v4) edge[double] node[below, font = \scriptsize] {$1_\mathcal{X}^*$} (v3);
\node (v5) at (-1.5,2.6) {};
\node (v6) at (-1.5,1.4) {};
\draw[->]  (v5) edge[double] node[right, font = \scriptsize] {$\alpha_f$} (v6);
\end{tikzpicture}
    \end{center}
    given, for $u \in Ff(x,y)$, $x \in \mathcal{X}_a$, $y \in \mathcal{X}_b$, by 
    $$\alpha_f(u) = (f, u)$$
\end{itemize}
By definition of $\mathcal{X}$ the commutatativity condition in Definition \ref{vtrans} holds. 

To verify the universal property of $(\mathcal{X}, \alpha)$ let $\mathcal{Y}$ be a category and $\beta : F \Rightarrow \mathcal{Y}$ be another vertical cone over $F$. We can define the functor $L : \mathcal{X} \rightarrow \mathcal{Y}$ simply by $L(a,x) = \beta_a(x)$ on objects and $L(f,u) = \beta_f(u)$ on morphisms. By construction $L \circ \alpha = \beta$ and hence the desired result.
\begin{flushright}
$\square$
\end{flushright}
\end{redi_proof}

\begin{remark}
The definition of $\mathcal{X}$ above is similar to that of the Grothendieck construction for diagrams of categories. As a matter of fact this is a generalization (see \cite{haderi2019higher} for a proof of the fact that the Grothendieck construction is a double colimit). 
\end{remark}

\begin{corollary}
Let $\mathcal{C}$, $\mathcal{D}$ be categories and $u$ be a $(\mathcal{C}, \mathcal{D})$-profunctor. When we regard $u$ as a functor $p: \textbf{col}(u) \rightarrow \Delta^1$, then in our notation the functor $\chi_p : \Delta^1 \rightarrow \textbf{Prof}_h$ is precisely the map that picks $u$. In this case we have
$$\textbf{dcolim}(\chi_p) \cong \textbf{col}(u)$$
\end{corollary}

\begin{corollary}
Let $\mathcal{A}$ be a category and $F: \mathcal{A} \rightarrow \textbf{Prof}_h$ be a lax diagram. Then $\textbf{dcolim}(F)$ is naturally equipped with a functor to $\mathcal{A}$ so that we have a functor
$$\textbf{dcolim} : [\mathcal{A}, \textbf{Prof}_h]_{lax} \rightarrow \textbf{Cat}/\mathcal{A}$$

\end{corollary}
\begin{redi_proof}
Let $\mathcal{X}$ be the category defined in Proposition \ref{model}. There is a functor $P : \mathcal{X} \rightarrow \mathcal{A}$ given by $P(a,x) = a$ on objects and $P(f,u) = f$ on morphisms. Clearly, the construction of $P$ is functorial. 
\begin{flushright}
$\square$
\end{flushright}
\end{redi_proof}

Now we can state and prove the theorem we aimed for.

\begin{theorem} \label{fun_class}
Let $\mathcal{A}$ be a category. The pair $(\chi, \textbf{dcolim})$
$$\chi : \textbf{Cat}/\mathcal{A} \leftrightarrows [\mathcal{A}, \textbf{Prof}_h]_{lax} : \textbf{dcolim}$$
as defined above is an equivalence of categories between the slice category of functors over $\mathcal{A}$ and the category of lax $\mathcal{A}$-indexed diagrams of profunctors with vertical transformations as morphisms. 

\end{theorem}

\begin{redi_proof}
Let $F : \mathcal{X} \rightarrow \mathcal{A}$ be a functor. By definition of $\chi_F$ and Proposition \ref{model} it is clear that
$$\textbf{dcolim} \chi_F \cong \mathcal{X}$$

On the other hand, consider a lax 2-functor $F: \mathcal{A} \rightarrow \textbf{Prof}_h$. Using the model $\mathcal{X}$ defined in Proposition \ref{model} for the double colimit of $F$, with the obvious functor $p: \mathcal{X} \rightarrow \mathcal{A}$, we actually have $\chi_p = F$.
\begin{flushright}
$\square$
\end{flushright}
\end{redi_proof}

In conclusion, functors into a category $\mathcal{A}$ are classified by unitally strict lax 2-functors $\mathcal{A} \rightarrow \textbf{Prof}_h$, and it is double category theory that provides the necessary organization, language and tools to understand this result. We will move on to the world of simplicial sets, and propose a way to understand their double category theory and state (and prove) an analog of Theorem \ref{fun_class}.

\begin{note}
It is reasonable to expect the equivalence $(\chi, \textbf{dcolim})$ in Theorem \ref{fun_class} to extend to an equivalence of 2-categories, given that $\textbf{Cat}/\mathcal{A}$ has a natural structure of a 2-category inherited from \textbf{Cat} (see also \cite{street2001powerful}). We discuss briefly and leave the details to the reader.

In general, for a double category $\mathbb{D}$, there is a natural 2-category structure on the category of horizontal $\mathcal{A}$-diagrams $[\mathcal{A}, \mathbb{D}_h]_{lax}$. For $F,G$ horizontal diagrams and $\alpha, \beta : F \Rightarrow G$ vertical transformations we can define horizontal modifications $\tau : \alpha \Rrightarrow \beta$ to consist of 2-cells $\tau_a$
\begin{center}
        \begin{tikzpicture}

\node (v1) at (-2.5,3) {$F(a)$};
\node (v4) at (-2.5,1) {$G(a)$};
\node (v2) at (-0.5,3) {$F(a)$};
\node (v3) at (-0.5,1) {$G(a)$};
\draw[-]  (v1) edge[double] (v2);
\draw[->]  (v2) edge node[right, font = \scriptsize] {$\beta_a$}  (v3);
\draw[->]  (v1) edge node[left, font = \scriptsize] {$\alpha_a$} (v4);
\draw[-]  (v4) edge[double]  (v3);
\node (v5) at (-1.5,2.6) {};
\node (v6) at (-1.5,1.4) {};
\draw[->]  (v5) edge[double] node[right, font = \scriptsize] {$\tau_a$} (v6);
\end{tikzpicture}
    \end{center}
for each $a \in \mathcal{A}$, which commute with the structure 2-cells of $\alpha$ and $\beta$. 
This way for $\mathbb{D} = \textbf{Prof}$ we obtain a 2-category $[\mathcal{A}, \textbf{Prof}_h]_{lax}$ which is equivalent to $\textbf{Cat}/\mathcal{A}$. 
\end{note}

\section{Higher correspondences and simplicial maps} \label{higher_corrs}

Our goal is to develop results similar to those we introduced above for simplicial sets. Making Lurie's prediction precise falls within this scope as well. In order to achieve this we have to develop double category theory for simplicial sets.

Of course, given that correspondences between simplicial sets cannot be tensored, we cannot have an actual double category whose objects are simplicial sets, vertical morphisms are simplicial maps and horizontal morphisms are correspondences. Nonetheless, we claim that a meaningful double category theory is achievable. 

The key idea we propose is to use simplicial categories instead of double categories. This will turn out to make sense from many angles:
\begin{itemize}
    \item Higher correspondences themselves look and feel like simplices. Once we define this simplicial structure it is possible to study functors valued in higher correspondences.
    \item Categories and functors are part of a 2-category, which we recorded as the vertical 2-category of the double category \textbf{Prof}.  Simplicial sets and simplicial maps are part of a simplicially enriched category. In general, we will see that a simplicial category, looked at from the double categorical perspective, has a vertical \textbf{sSet}-category as part of it. Thus the latter become good candidates for our purposes.
    \item The theory of double colimits is available in simplicial categories.
    \item Some of the expected results can be proved.
\end{itemize}

\subsection{A simplicial category for higher correspondences} \label{higher}

\begin{definition}
An $n$-correspondence of simplicial sets is a simplicial set $X$ equipped with a map $p: X \rightarrow \Delta^n$.
\end{definition}

As mentioned in the introduction, the structure map $p$ allows us to think of $X$ as a collage of simplicial sets. First, we can think of the fibers $X_i = p^{-1}(i)$ as the vertices of $X$. Then we can describe $X$ to be formed from the $X_i$'s by adding new simplices in the increasing direction (Figure \ref{2-corr}). This picture indicates that an $n$-correspondence itself should be regarded as some sort of $n$-simplex.

\begin{figure}
    \centering
     \begin{tikzpicture}
\draw  (-2.4,2.8) ellipse (0.4 and 0.6);
\draw  (0,4.4) ellipse (0.4 and 0.6);
\draw  (2.4,2.8) ellipse (0.4 and 0.6);

\node (v1) at (-2.4,3.2) {};
\node (v4) at (-2.4,2.4) {};
\node (v2) at (0,4.8) {};
\node (v5) at (0,4) {};
\node (v3) at (2.4,3.2) {};
\node (v6) at (2.4,2.4) {};
\draw[->]  (v1) edge (v2);
\draw[->]  (v2) edge (v3);
\draw[->]  (v1) edge (v3);
\draw[->]  (v4) edge (v5);
\draw[->]  (v5) edge (v6);
\draw[->]  (v4) edge (v6);
\node at (-3.2,3.4) {$X_0$};
\node at (0.8,4.8) {$X_1$};
\node at (3.2,3.4) {$X_2$};
\end{tikzpicture}
    \caption{A 2-correspondence}
    \label{2-corr}
\end{figure}

Indeed we may define a simplicial category $\textbf{sSet}^\sharp$, meaning a category-valued mapping
$$\textbf{sSet}^\sharp : \Delta^{op} \rightarrow \textbf{Cat}$$
with category of $n$-simplices being the category of $n$-correspondences, i.e. 
$$\textbf{sSet}^\sharp_n = \textbf{sSet}/\Delta^n$$
The action of $\Delta$ is given by pullbacks. Given a morphism $\theta : [m] \rightarrow [n]$ in $\Delta$ there is a functor $\theta^* : \textbf{sSet}^\sharp_n \rightarrow \textbf{sSet}^\sharp_m$ which maps an $n$-correspondence $p: X \rightarrow \Delta^n$  to the $m$-correspondence $\theta^*p : \theta^*X \rightarrow \Delta^m$ obtained pulling back along $\theta$

\begin{center}
	\begin{tikzpicture}
	
	\node (v4) at (-3.5,4) {$\theta^*X$};
	\node (v1) at (-2,4) {$X$};
	\node (v3) at (-3.5,2.5) {$\Delta^{m}$};
	\node (v2) at (-2,2.5) {$\Delta^n$};
	\draw[->]  (v1) edge node[right, font= \scriptsize]{$p$} (v2);
	\draw[->]  (v3) edge node[above, font= \scriptsize]{$\theta$} (v2);
	\draw[->, dashed]  (v4) edge node[left, font= \scriptsize]{$\theta^*p$} (v3);
	\draw[->, dashed]  (v4) edge (v1);
	\end{tikzpicture}
\end{center}

The pasting laws for pullback squares assure that this simplicial structure is well defined, but only up to natural isomorphism, i.e. $\textbf{sSet}^\sharp$ is a weak 2-functor. It is only fair that when we study simplicial objects in a 2-category we are faced with weak functoriality. Lets illustrate this structure (and its weakness) by computing $\theta^*$ for some specific examples of $\theta$.

\begin{example}
Let $p : X \rightarrow \Delta^n$ be an $n$-correspondence and let $d^i : [n-1] \rightarrow [n]$ be the $i$-th face inclusion. Then the face $d_i^*X$ is simply the simplicial set obtained by ``deleting" the vertex $X_i$. In other words $d_i^*X$ is the largest simplicial subset of $X$ not containing the $0$-simplices of $X_i$.
\end{example}

\begin{example} \label{degeneracy}
Let $X$ be a simplicial set and let $t^n : [n] \rightarrow [0]$ be the terminal map in $\Delta$. Then, regarding $X$ as a $0$-correspondence, the $n$-correspondence $t_n^*X$ is simply the cylinder $X \times \Delta^n$ with structure map the projection to $\Delta^n$. 

This fact is consistent with the theory of profunctors where, for a category $\mathcal{C}$, the identity profunctor, which is interpreted as a degeneracy in the language of simplices, is represented by the cylinder $\mathcal{C} \times \Delta^1$. 
\end{example}

\begin{remark}
Example \ref{degeneracy} serves to illustrate the weakness in $\textbf{sSet}^\sharp$. When we consider the cylinder $X \times \Delta^1$ for a simplicial set $X$, which is the degeneracy $s_0^*X$, its top and bottom, which are the faces $d_0^*s_0^*X$ and $d_1^*s_0^*X$, are copies of $X$ but not identically $X$. 
\end{remark}

\begin{example}
In general, for an $n$-correspondence $X \rightarrow \Delta^n$ the degeneracy $s_i^*X$ may be thought of as obtained by ``extending" $X$ with a cylinder along $X_i$. For clarity we will work out an example. Let $X_0 = (a)$ be a copy of $\Delta^0$ and $X_1 = (b \rightarrow c)$ be a copy of $\Delta^1$. Consider $X = \Delta^2$ as a 1-correspondence between them via the map $s_1 : \Delta^2 \rightarrow \Delta^1$
\begin{center}
	\begin{tikzpicture}
	
	\draw  (-4.1,3.5) ellipse (0.3 and 0.5);
	\draw  (-2.5,3.5) ellipse (0.5 and 1);
	\node (v2) at (-2.5,4) {$b$};
	\node (v3) at (-2.5,3) {$c$};
	\node (v1) at (-4.1705,3.4784) {$a$};
	\draw[->]  (v1) edge (v2);
	\draw[->]  (v1) edge (v3);
	\draw[->]  (v2) edge (v3);
	\node at (-4.1165,2.7043) {$X_0$};
	\node at (-2.4963,2.1462) {$X_1$};
	\end{tikzpicture}
\end{center}
Then $s_0^*X$ and $s_1^*X$ are depicted as
\begin{center}
	
	\begin{tikzpicture}
	
	\draw  (-4.5,3.5) ellipse (0.3 and 0.5);
	\draw  (-2,3.5) ellipse (0.5 and 1);
	\node (v2) at (-2,4) {$b$};
	\node (v3) at (-2,3) {$c$};
	\node (v1) at (-4.6,3.5) {$a_0$};
	\draw[->]  (v1) edge (v2);
	\draw[->]  (v1) edge (v3);
	\draw[->]  (v2) edge (v3);

	\draw  (-3.2,5.5) ellipse (0.3 and 0.5);
	\node (v4) at (-3.3,5.5) {$a_1$};
	\draw[->]  (v1) edge (v4);
	\draw[->]  (v4) edge (v2);
	\draw[->]  (v4) edge (v3);
	\end{tikzpicture}\ \ \ \ \ \ \ \
	\begin{tikzpicture}
	
	\draw  (-4.1,3.5) ellipse (0.3 and 0.5);
	\draw  (-2.4,5.2) ellipse (0.5 and 1);
	\node (v2) at (-2.4,5.6) {$b_0$};
	\node (v3) at (-2.4,4.6) {$c_0$};
	\node (v1) at (-4.1705,3.4784) {$a$};
	\draw[->]  (v1) edge (v2);
	\draw[->]  (v1) edge (v3);
	\draw[->] (v2) edge (v3);
	
	\draw  (-1,3) ellipse (0.5 and 1);
	\node (v4) at (-1,3.5) {$b_1$};
	\node (v5) at (-1,2.5) {$c_1$};
	\draw[->]  (v2) edge (v4);
	\draw[->]  (v4) edge (v5);
	\draw[->]  (v3) edge (v5);
	\draw[->]  (v2) edge (v5);
	\draw[->]  (v1) edge (v4);
	\draw[->]  (v1) edge (v5);
	\end{tikzpicture}
\end{center}
\end{example}

Let $\infty\textbf{Cat}$ be the category of $\infty$-categories (meaning the full subcategory of \textbf{sSet} consisting of quasi-categories). We can define the simplicial category $\infty\textbf{Cat}^\sharp$ exactly as above with 
$$\infty\textbf{Cat}^\sharp_n = \infty\textbf{Cat}/\Delta^n$$
The above prescription works well for $\infty$-categories because of the following easy-to-prove two facts (\cite{lurie2009higher})
\begin{itemize}
    \item[i)] a map of simplicial sets $p: X \rightarrow A$, where $A$ is a category, is an inner fibration if and only if $X$ is an $\infty$-category
    \item[ii)] inner fibrations are stable under products and pullbacks 
\end{itemize}

In a similar fashion we could define $\textbf{Cat}^\sharp$ or $\textbf{Top}^\sharp$, where \textbf{Cat} stands for categories and \textbf{Top} for topological spaces. Other interesting examples of simplicial categories not involving maps to $\Delta^n$ can be defined. 
For instance, we can take the horizontal lax nerve (also known as Duskin nerve in 2-category theory \cite{duskin2002simplicial}) of any double category $\mathbb{D}$ to obtain a simplicial category $N_h(\mathbb{D})$ with category of $n$-simplices $[\Delta^n, \mathbb{D}_h]_{lax}$. In light of Theorem \ref{fun_class} we have $\textbf{Cat}^\sharp \cong N_h(\textbf{Prof})$.

\subsection{The analogy with double categories}
\label{analogy}

We will briefly discuss the double categorical aspects of simplicial categories. Let 
$$\mathbb{E} : \Delta^{op} \rightarrow \textbf{Cat}$$
be a simplicial category. Given our examples from the previous section we will always assume $\mathbb{E}$ to be weak by default. Composing $\mathbb{E}$ with the nerve functor $N : \textbf{Cat} \rightarrow \textbf{sSet}$ allows us to visualise $\mathbb{E}$ as a bisimplicial set, consisting of simplices in both directions tied to each other by bisimplices. One has to keep in mind though that $N \circ \mathbb{E}$ is a weak functor from $\Delta^{op}$ into a $2$-category of simplicial sets, while bisimplicial sets are strict mappings.

More explicitly $\mathbb{E}$ may be depicted as a 2-fold categorical structure with:
\begin{itemize}
    \item objects those of $\mathbb{E}_0$
    \item vertical morphisms those of $\mathbb{E}_0$
    \item the objects of $\mathbb{E}_n$ for various $n$ as horizontal simplices
    \item morphisms of $\mathbb{E}_n$ as prism-shaped cells  of various dimensions. For example a morphism $f: x \rightarrow y$ in $\mathbb{E}_1$ looks like a square
    \begin{center}
	\begin{tikzpicture}
	
	\node (v1) at (-3.5,3.5) {$x_0$};
	\node (v3) at (-2,3.5) {$x_1$};
	\node (v2) at (-3.5,2) {$y_0$};
	\node (v4) at (-2,2) {$y_1$};
	\draw[->]  (v1) edge node[left, font= \scriptsize]{$f_0$} (v2);
	\draw[->]  (v1) edge node[above, font= \scriptsize]{$x$} (v3);
	\draw[->]  (v2) edge node[below, font= \scriptsize]{$y$} (v4);
	\draw[->]  (v3) edge node[right, font= \scriptsize]{$f_1$} (v4);
	
	\node at (-2.7505,2.7145) {$f$};
	\end{tikzpicture}
\end{center}
with $f_0$ and $f_1$ being the faces of $f$. If $f \in \mathbb{E}_2$ it looks like a prism 
\begin{center}
	\begin{tikzpicture}

	\node (v1) at (-4,2.5) {$x_0$};
	\node (v2) at (-2.5,3.5) {$x_1$};
	\node (v3) at (-1,2.5) {$x_2$};
	\node (v4) at (-4,0.5) {$y_0$};
	\node (v5) at (-2.5,1.5) {$y_1$};
	\node (v6) at (-1,0.5) {$y_2$};
	\draw[->]  (v1) edge (v2);
	\draw[->]  (v2) edge (v3);
	\draw[->]  (v1) edge (v3);
	\draw[->]  (v1) edge node[left, font= \scriptsize]{$f_0$} (v4);
	\draw[->]  (v2) edge node[left, font= \scriptsize]{$f_1$} (v5);
	\draw[->]  (v3) edge node[right, font= \scriptsize]{$f_2$} (v6);
	\draw[->]  (v4) edge (v5);
	\draw[->]  (v5) edge (v6);
	\draw[->]  (v4) edge (v6);
	\end{tikzpicture}
\end{center}
And so on in higher dimensions.
\end{itemize}

If $\mathbb{D}$ is a double category we extract its vertical 2-category $\mathbb{D}_v$ by taking cells whose top and bottom horizontal morphisms are identities. In a simplicial context the role of identities is played by degeneracies. If $\mathbb{E}$ is a simplicial category, when we consider only cells in which the top and bottom are totally degenerate we are left with a simplicially enriched category $\mathbb{E}_v$.

More precisely, let $t^n : [n] \rightarrow [0]$ be the terminal map in $\Delta$. 
Then $\mathbb{E}_v$ is the \textbf{sSet}-category whose objects are those of $\mathbb{E}$, and whose $n$-simplices of the mapping space $\mathbb{E}_v(x,y)$, for $x,y \in \mathbb{E}$, are given by
$$\mathbb{E}_v(x,y)_n = \mathbb{E}_n(t_n^*x, t_n^*y)$$
Faces and degeneracies are defined in the obvious way using the simplicial structure of $\mathbb{E}$. 

What the double categorical lens unveiled for us is that 
\begin{quote}
    ``Simplicial categories are a 2-fold variant of \textbf{sSet}-categories"
\end{quote}
Categories, functors and natural transformations are part of a 2-category which serves as the vertical part of \textbf{Prof}. Simplicial sets, simplicial maps and homotopies of various dimensions are part of a \textbf{sSet}-category \underline{\textbf{sSet}}. Recall that for two simplicial sets $X$ and $Y$ the set of $n$-simplices of their mapping space is defined to be
$$\underline{\textbf{sSet}}(X,Y)_n = \textbf{sSet}(X \times \Delta^n, Y)$$

\underline{\textbf{sSet}} is the vertical \textbf{sSet}-category of $\textbf{sSet}^\sharp$. Indeed, let $X$ and $Y$ be simplicial sets. Then we have 
\begin{align*}
    \textbf{sSet}^\sharp_v(X,Y)_n &= \textbf{sSet}^\sharp_n(t_n^*X, t_n^*Y) \\
    &= \textbf{sSet}/\Delta^n(X \times \Delta^n, Y \times \Delta^n) \\
    &\cong \textbf{sSet}(X \times \Delta^n, Y)
\end{align*}
We see that under our observations $\textbf{sSet}^\sharp$ is indeed analogous to \textbf{Prof}. 

The analogy goes further. With degeneracies playing the role of identities and simplices the role of composites we can easily extend the theory of double colimits as well to simplicial categories. 
More precisely, if $\mathbb{E}$ is a simplicial category and $A$ is a simplicial set we can study double colimits of horizontal $A$-indexed diagrams in $\mathbb{E}$. 

When we regard $A$ as a discrete simplicial category then a horizontal $A$-diagram in $\mathbb{E}$ can be thought of as a weak natural transformation (between mappings $\Delta^{op} \rightarrow \textbf{Cat}$) $F: A \rightarrow \mathbb{E}$. This makes sense because $F$ assigns to each simplex $\sigma \in A_n$ an object $F(\sigma) \in \mathbb{E}_n$ and from our viewpoint the latter is an $n$-simplex in the horizontal direction of $\mathbb{E}$. $F$ being weak simply says that for morphisms $\theta : [m] \rightarrow [n]$ in $\Delta$ we have isomorphisms $F(\sigma\theta) \cong \mathbb{E}(\theta)F(\sigma)$ in the category $\mathbb{E}_m$ (rather than strict equalities) which satisfy some coherence laws. Again, besides being appropriate for studying examples, this makes sense conceptually as well because we take into account the fact that $\mathbb{E}_n$'s are categories rather than sets. 

In this vein, given two horizontal diagrams $F, G: A \rightarrow \mathbb{E}$, a vertical transformation $\alpha: F \Rightarrow G$ between them is simply a modification between weak natural transformations. $\alpha$ assigns to each $\sigma \in A_n$ a morphism $\alpha_\sigma : F(\sigma) \rightarrow G(\sigma)$ in $\mathbb{E}_n$, the latter being thought of as cells in the vertical direction of $\mathbb{E}$. Moreover, the structure maps $\alpha_\sigma$ have to respect the simplicial structure, i.e. for all $\theta : [m] \rightarrow [n]$ in $\Delta$ we have $\alpha_{\sigma\theta} = \mathbb{E}(\theta)\alpha_\sigma$. Hence, horizontal $A$-indexed diagrams in $\mathbb{E}$ form a category which we denote $[A, \mathbb{E}]$. 

\begin{lemma} \label{yoneda}
Let $\mathbb{E}$ be a simplicial category. Then there is an equivalence of categories
$$[\Delta^n, \mathbb{E}] \cong \mathbb{E}_n$$
\end{lemma}
\begin{redi_proof}
This is a direct implementation of the Yoneda Lemma for weak 2-categories (see Lemma 8.3.11 in \cite{johnson20212}).
\begin{flushright}
$\square$
\end{flushright}
\end{redi_proof}

The statement of Lemma \ref{yoneda} is fairly intuitive and it tells us that objects in $\mathbb{E}_n$ can be regarded up to isomorphism as horizontal diagrams $\Delta^n \rightarrow \mathbb{E}$. We will use this identification freely.

Now we have the tools to define the double colimit of a horizontal diagram $F: A \rightarrow \mathbb{E}$. For an object $x \in \mathbb{E}_0$ there is a constant $A$-diagram 
$$c_x : A \rightarrow \Delta^0 \xrightarrow{x} \mathbb{E}$$
given by composing the terminal map from $A$ to $\Delta^0$ with the map from $\Delta^0$ to $\mathbb{E}$ which picks $x$. Then we define the category of vertical cones over $F$, denoted $\textbf{vCone}(F)$, with
\begin{itemize}
    \item object pairs $(x, \alpha)$ where $x \in \mathbb{E}_0$ and $\alpha: F \Rightarrow c_x$ is a vertical transformation into the constant diagram 
    \item morphisms $f: (x, \alpha) \rightarrow (y, \beta)$ between two objects being morphisms $f: x \rightarrow y$ in $\mathbb{E}_0$ such that $c_f \circ \alpha = \beta$ in $[A, \mathbb{E}]$, where $c_f : c_x \Rightarrow c_y$ is the vertical transformation between constant diagrams induced by $f$
\end{itemize}

Again observe that just like for double categories we ask of morphisms in $\textbf{vCone}(F)$ to come from $\mathbb{E}_0$ and not simply be morphisms between constant diagrams commuting with the structure maps. 

\begin{definition}
Let $A$ be a simplicial set, $\mathbb{E}$ a simplicial category and $F : A \rightarrow \mathbb{E}$ a horizontal diagram. The double colimit of $F$ is an object $\textbf{dcolim}F \in \mathbb{E}_0$ equipped with a vertical transformation $F \Rightarrow c_{\textbf{dcolim}F}$ which is initial in $\textbf{vCone}(F)$.
\end{definition}

Alternatively, we can define the double colimit simultaneously for all diagrams in $[A, \mathbb{E}]$ by asking for a functor \textbf{dcolim} which is left adjoint to the constant diagram functor $c: \mathbb{E}_0 \rightarrow [A, \mathbb{E}]$. If each diagram has a double colimit then these definitions are equivalent.

\begin{note}
Double colimits in simplicial categories as we have defined them above should extend those in double categories in the sense that they should be preserved by horizontal lax nerves. However, this is too much of a formality to be proven here since we focus on the effectiveness of double colimits rather than their formal properties.
\end{note}

\subsection{Double colimits of diagrams of correspondences}

Now we prove an analog of Proposition \ref{model} for simplicial categories, thus providing an explicit model for all double colimits in $\textbf{sSet}^\sharp$. 

\begin{construction} \label{X}
Let $F: A \rightarrow \textbf{sSet}^\sharp$ be a diagram of correspondences. For a simplex $\sigma \in A_n$ we denote its image under $F$ by $p_\sigma : F(\sigma) \rightarrow \Delta^n$. 

We construct a simplicial set $X$ with set of $n$-simplices defined to be
$$X_n = \{(\sigma, x) \ | \ \sigma \in A_n, x \in F(\sigma)_n  \ \text{section of} \ p_\sigma \}$$
The action of $\Delta$ is given as follows. Let $\theta : [m] \rightarrow [n]$ be a morphism in $\Delta$ and $(\sigma, x) \in X_n$. By definition of $\textbf{sSet}^\sharp$ we have a pullback square

\begin{center}
    \begin{tikzpicture}

\node (v1) at (-2,2) {$F(\sigma \theta)$};
\node (v2) at (0,2) {$F(\sigma)$};
\node (v3) at (-2,0) {$\Delta^m$};
\node (v4) at (0,0) {$\Delta^n$};
\node (v5) at (-3,3) {$\Delta^m$};
\draw[->]  (v1) edge (v2);
\draw[->]  (v1) edge node[right, font= \scriptsize]{$p_{\sigma \theta}$} (v3);
\draw[->]  (v2) edge node[right, font= \scriptsize]{$p_\sigma$} (v4);
\draw[->]  (v3) edge node[below, font= \scriptsize]{$\theta$} (v4);
\draw[->, bend left = 15]  (v5) edge node[above, font= \scriptsize]{$x\theta$} (v2);
\draw[->, bend right = 15]  (v5) edge node[left, font= \scriptsize]{$1_{\Delta^m}$} (v3);
\draw[->, dashed]  (v5) edge node[right, font= \scriptsize]{$\theta^* x$} (v1);
\end{tikzpicture}
\end{center}
Note that we are taking advantage of the isomorphism $F(\sigma \theta) \cong \theta^*F(\sigma)$ in $\textbf{sSet}/\Delta^m$ in asserting that the square on the right is a pullback square. We will tacitly keep using the structure isomorphisms of $F$ throughout this construction (and the rest of this work) in order not to overload our diagrams.

Then we obtain an induced map $\theta^*x : \Delta^m \rightarrow F(\sigma \theta)$ with components $x\theta$ and $1_{\Delta^m}$, and hence a section of $p_{\sigma \theta}$. Thus we define 
$$(\sigma, x)\theta = (\sigma\theta, \theta^*x)$$
This construction is clearly functorial by pasting properties of pullback squares.

Next we construct a vertical transformation $\alpha : F \Rightarrow c_X$, thus exhibiting $X$ as a vertical cone over $F$.
Such an $\alpha$ will consist of components $\alpha_\sigma : F(\sigma) \rightarrow t_n^*X$ in $\textbf{sSet}/\Delta^n$ for $\sigma \in A_n$, $[n] \in \Delta$. Recall from Example \ref{degeneracy} that $t_n^*X$ is simply the $n$-correspondence given by the projection $X \times \Delta^n \rightarrow \Delta^n$. Since we have a natural isomorphism 
$$\textbf{sSet}/\Delta^n(F(\sigma), X \times \Delta^n) \cong \textbf{sSet}(F(\sigma), X)$$
it is enough to define maps of simplicial sets $\alpha_\sigma : F(\sigma) \rightarrow X$. 

So let $\sigma \in A_n$ and $y : \Delta^k \rightarrow F(\sigma)$ be a $k$-simplex. Pulling back along $p_\sigma y $ 

\begin{center}
    \begin{tikzpicture}

\node (v1) at (-2,2) {$F(\sigma p_\sigma y)$};
\node (v2) at (0,2) {$F(\sigma)$};
\node (v3) at (-2,0) {$\Delta^k$};
\node (v4) at (0,0) {$\Delta^n$};
\node (v5) at (-3,3) {$\Delta^k$};
\draw[->]  (v1) edge (v2);
\draw[->]  (v1) edge node[right, font= \scriptsize]{$p_{\sigma p_\sigma y}$} (v3);
\draw[->]  (v2) edge node[right, font= \scriptsize]{$p_\sigma$} (v4);
\draw[->]  (v3) edge node[below, font= \scriptsize]{$p_\sigma y$} (v4);
\draw[->, bend left = 15]  (v5) edge node[above, font= \scriptsize]{$y$} (v2);
\draw[->, bend right = 15]  (v5) edge node[left, font= \scriptsize]{$1_{\Delta^k}$} (v3);
\draw[->, dashed]  (v5) edge node[right, font= \scriptsize]{$\Tilde{y}$} (v1);
\end{tikzpicture}
\end{center}

we obtain an induced map $\Tilde{y} : \Delta^k \rightarrow F(\sigma p_\sigma y)$ with components $1_{\Delta^k}$ and $y$. By construction $\Tilde{y}$ is a section of $p_{\sigma p_\sigma y}$ and hence defines a $k$-simplex in $X$. Thus we define 
$$\alpha_{\sigma,k}(y) = (\sigma p_\sigma y, \Tilde{y})$$

We first check that $\alpha_\sigma$ is a well-defined map of simplicial sets. So let $\tau : [l] \rightarrow [k]$ be a morphism in $\Delta$. We have the following commutative diagram 

\begin{center}
    \begin{tikzpicture}

\node (v1) at (-3,2) {$F(\sigma p_\sigma y \tau)$};
\node (v4) at (-3,0) {$\Delta^l$};
\node (v2) at (-1,2) {$F(\sigma p_\sigma y)$};
\node (v5) at (-1,0) {$\Delta^k$};
\node (v3) at (1,2) {$F(\sigma)$};
\node (v6) at (1,0) {$\Delta^n$};
\node (v7) at (-4,3) {$\Delta^l$};
\node (v8) at (-2,3) {$\Delta^k$};
\draw[->]  (v1) edge (v2);
\draw[->]  (v2) edge (v3);
\draw[->]  (v1) edge node[right, font= \scriptsize]{$p_{\sigma p_\sigma y \tau}$} (v4);
\draw[->]  (v2) edge node[right, font= \scriptsize]{$p_{\sigma p_\sigma y}$} (v5);
\draw[->]  (v3) edge node[right, font= \scriptsize]{$p_\sigma$} (v6);
\draw[->]  (v4) edge node[below, font= \scriptsize]{$\tau$} (v5);
\draw[->]  (v5) edge node[below, font= \scriptsize]{$p_\sigma y$} (v6);
\draw[->]  (v7) edge node[above, font= \scriptsize]{$\tau$} (v8);
\draw[->, bend left = 15]  (v8) edge node[above, font= \scriptsize]{$y$} (v3);
\draw[->, dashed]  (v8) edge node[right, font= \scriptsize]{$\Tilde{y}$} (v2);
\draw[->, dashed]  (v7) edge node[right, font= \scriptsize]{$\tau^* \Tilde{y}$} (v1);
\draw[->, bend right= 15]  (v7) edge node[left, font= \scriptsize]{$1_{\Delta_l}$} (v4);
\end{tikzpicture}
\end{center}

The two squares on the bottom are pullback squares and the dashed arrows represent maps induced by the universal property of pullbacks. Because pasting pullback squares results in a pullback square we have the following equality in $X_l$
$$(X(\tau) \circ \alpha_{\sigma,k})(y) = (\sigma p_\sigma y \tau, \tau^* \Tilde{y}) = (\alpha_{\sigma, l} \circ F(\sigma)(\tau))(y)$$
Hence $\alpha_\sigma$ is a well-defined map of simplicial sets. 

Next we check that $\alpha_\sigma$ is well-defined with respect to $\sigma$. So let $\theta : [m] \rightarrow [m]$ be a morphism in $\Delta$. 
Let $f_\sigma : F(\sigma \theta) \rightarrow F(\sigma)$ be the top map in the pullback square

\begin{center}
    \begin{tikzpicture}

\node (v1) at (-2,2) {$F(\sigma \theta)$};
\node (v2) at (0,2) {$F(\sigma)$};
\node (v3) at (-2,0) {$\Delta^m$};
\node (v4) at (0,0) {$\Delta^n$};

\draw[->]  (v1) edge node[above, font= \scriptsize]{$f_\sigma$} (v2);
\draw[->]  (v1) edge node[left, font= \scriptsize]{$p_{\sigma \theta}$} (v3);
\draw[->]  (v2) edge node[right, font= \scriptsize]{$p_\sigma$} (v4);
\draw[->]  (v3) edge node[below, font= \scriptsize]{$\theta$} (v4);

\end{tikzpicture}
\end{center}

Note that it is enough to prove $\alpha_\sigma f_\sigma = \alpha_{\sigma \theta}$ to show that $\alpha$ is well-defined. So let $y : \Delta^k \rightarrow F(\sigma \theta)$ be a $k$-simplex. We have the following diagram

\begin{center}
    \begin{tikzpicture}

\node (v1) at (-3,2) {$F(\sigma \theta p_{\sigma\theta} y)$};
\node (v4) at (-3,0) {$\Delta^k$};
\node (v2) at (-1,2) {$F(\sigma \theta)$};
\node (v3) at (1,2) {$F(\sigma)$};
\node (v5) at (-1,0) {$\Delta^m$};
\node (v6) at (1,0) {$\Delta^m$};
\node (v7) at (-4,3) {$\Delta^k$};
\draw[->]  (v1) edge (v2);
\draw[->]  (v2) edge node[above, font= \scriptsize]{$f_\sigma$} (v3);
\draw[->]  (v1) edge node[right, font= \scriptsize]{$p_{\sigma \theta p_{\sigma\theta} y}$} (v4);
\draw[->]  (v2) edge node[right, font= \scriptsize]{$p_{\sigma \theta}$} (v5);
\draw[->]  (v3) edge node[right, font= \scriptsize]{$p_\sigma$} (v6);
\draw[->]  (v4) edge node[below, font= \scriptsize]{$p_{\sigma \theta} y$} (v5);
\draw[->]  (v5) edge node[below, font= \scriptsize]{$\theta$} (v6);
\draw[->, bend left = 15]  (v7) edge node[above, font= \scriptsize]{$f_\sigma y$} (v3);
\draw[->, bend left = 10]  (v7) edge node[right, font= \scriptsize]{$y$} (v2);
\draw[->, dashed]  (v7) edge node[right, font= \scriptsize]{$\Tilde{y}$} (v1);
\draw[->, bend right = 15]  (v7) edge node[left, font= \scriptsize]{$1_{\Delta^k}$} (v4);
\end{tikzpicture}
\end{center}

where the bottom squares are pullback squares and $\Tilde{y}$ is the induced map. By definition of $\alpha$ and the pasting property of pullback squares we have 
$$\alpha_{\sigma \theta} (y) = (\sigma \theta p_{\sigma \theta} y, \Tilde{y}) = \alpha_\sigma (f_\sigma y)$$
Thus $\alpha$ is a well-defined vertical transformation $F \Rightarrow c_X$.

\end{construction}

\begin{proposition} \label{smodel}
Let $F : A \rightarrow \textbf{sSet}^\sharp$ be a horizontal diagram of correspondences. Then the pair $(X,\alpha)$ constructed above is a double colimit of $F$.
\end{proposition}

\begin{redi_proof}
We simply check the universality of $(X, \alpha)$. Let $(Y, \beta) \in \textbf{vCone}(F)$ be another such pair. We can define a map of simplicial sets $f : X \rightarrow Y$ by 
$$f(\sigma,x) = \beta_\sigma x$$
for all $(\sigma, x) \in X_n$. By definition $c_f \circ \alpha = \beta$.
\begin{flushright}
$\square$
\end{flushright}
\end{redi_proof}

\begin{remark}
While double colimits in \textbf{Prof} are a generalization of the Grothendieck construction for categories, it is reasonable to expect double colimits in $\textbf{sSet}^\sharp$ to be a generalization of homotopy colimits (in the sense of simplicially enriched categories). See \cite{haderi2019higher} for a (not yet complete) presentation of this idea.  
\end{remark}

\begin{corollary}
Let $F : \Delta^n \rightarrow \textbf{sSet}^\sharp$ be the map which picks the $n$-correspondence $p : X \rightarrow \Delta^n$. Then 
$$\textbf{dcolim}(F) \cong X$$
\end{corollary}

\begin{corollary} \label{slice}
For a simplicial set $A$, the double colimit of a diagram $F: A \rightarrow \textbf{sSet}^\sharp$ is naturally equipped with a map to $A$ so that we have a functor 
$$\textbf{dcolim} : [A, \textbf{sSet}^\sharp] \rightarrow \textbf{sSet}/A$$
\end{corollary}
\begin{redi_proof}
$\textbf{dcolim}(F) \cong X$, where $X$ is the simplicial set defined in Construction \ref{X}. Then there is a ``forgetful" map $f: X \rightarrow A$ given by $f(\sigma, x) = \sigma$. Functionality follows easily. 
\begin{flushright}
$\square$
\end{flushright}
\end{redi_proof}

There is a more conceptual way of defining and understanding $X$ as being obtained via a bundle. We discuss briefly.

\begin{definition}
A pointed $n$-correspondence is an $n$-correspondence $p: X \rightarrow \Delta^n$ equipped with a choice of section $x: \Delta^n \rightarrow X$.
\end{definition}

Note that a pointed correspondence can be naturally interpreted as a correspondence of pointed simplicial sets as well. Given a pointed $n$-correspondence $(X,x)$, the fibers $X_i = p^{-1}(i)$ will be simplicial sets with a chosen basepoint. Then we can think of $X$ as a simplicial set formed by adding new simplices joining vertices in the $X_i$'s in the increasing direction in such a way that we form a distinguished $n$-simplex $x$ joining the basepoints. 

There is a simplicial category $\textbf{sSet}^\sharp_*$ whose category of $n$-simplices is the category of pointed $n$-correspondences, defined exactly as $\textbf{sSet}^\sharp$. Here, a morphism of pointed correspondences is simply a simplicial map which preserves the choice of section. We have a ``forgetful map"
$$U : \textbf{sSet}^\sharp_* \rightarrow \textbf{sSet}^\sharp$$

Given a horizontal diagram $F: A \rightarrow \textbf{sSet}^\sharp$ the map $f: X \rightarrow A$ defined above fits in a commutative square

\begin{center}
    \begin{center}
    \begin{tikzpicture}

\node (v1) at (-2.5,3) {$X$};
\node (v4) at (-2.5,1) {$A$};
\node (v2) at (-0.5,3) {$\textbf{sSet}^\sharp_*$};
\node (v3) at (-0.5,1) {$\textbf{sSet}^\sharp$};
\draw[->]  (v1) edge node[above, font = \scriptsize]{$F_*$} (v2);
\draw[->]  (v2) edge node[right, font = \scriptsize]{$U$} (v3);
\draw[->]  (v1) edge node[left, font = \scriptsize]{$f$} (v4);
\draw[->]  (v4) edge node[below, font = \scriptsize]{$F$} (v3);
\end{tikzpicture}
\end{center}
\end{center}

The top map is obviously defined by $F_*(\sigma , x) = (F(\sigma), x)$. Moreover, it is not difficult to verify that the above is a pullback square of simplicial categories. 

\begin{remark}
Lurie points out in \cite[p. ~87]{lurie2009higher} that a higher category of correspondences cannot be realized as an $\infty$-category because higher morphisms do not have to be invertible. We regard $\textbf{sSet}^\sharp$ as a model for this higher structure, at least from the general point of view on the subject presented in \cite{leinster2004higher}, where higher categories are conceived of as structures comprised of objects, arrows and higher arrows of some sort. In a simplicial category $\mathbb{E}$ the morphisms in $\mathbb{E}_n$ play the role of higher arrows. 
\end{remark}

\subsection{Classifying simplicial maps and inner fibrations} \label{class_results}

Now that we have enough language at our disposal we can state and prove our results. As it happens a lot in category theory, once the correct setup is in place the proofs are really easy. 

Let $f: X \rightarrow A$ be a map of simplicial sets. For a simplex $\sigma \in A_n$ let $X_\sigma = f^{-1}(\sigma)$ be the preimage, given by the pullback square
\begin{center}
    \begin{tikzpicture}

\node (v1) at (-2.5,3) {$X_\sigma$};
\node (v4) at (-2.5,1) {$\Delta^n$};
\node (v2) at (-0.5,3) {$X$};
\node (v3) at (-0.5,1) {$A$};
\draw[->, dashed]  (v1) edge (v2);
\draw[->]  (v2) edge node[right, font = \scriptsize]{$f$} (v3);
\draw[->, dashed]  (v1) edge  (v4);
\draw[->]  (v4) edge node[above, font = \scriptsize]{$\sigma$} (v3);
\end{tikzpicture}
\end{center}
By definition $X_\sigma$ is an $n$-correspondence. The assignment $\sigma \mapsto X_\sigma$ produces a horizontal diagram of correspondences
$$\chi_f : A \rightarrow \textbf{sSet}^\sharp$$
Moreover, if we have a commutative triangle 

\begin{center}
    \begin{tikzpicture}

\node (v2) at (-1,1) {$A$};
\node (v1) at (-2,2) {$X$};
\node (v3) at (0,2) {$Y$};
\draw[->]  (v1) edge node[left, font = \scriptsize]{$f$} (v2);
\draw[->]  (v3) edge node[right, font = \scriptsize]{$g$} (v2);
\draw[->]  (v1) edge node[above, font = \scriptsize]{$h$} (v3);
\end{tikzpicture}
\end{center}
we obtain a vertical transformation $ \chi_f \Rightarrow \chi_g$ with components the maps $X_\sigma \rightarrow Y_\sigma$ induced by pulling back along simplices $\sigma$ in $A$. Thus the assignment $f \mapsto \chi_f$ is also functorial in the sense that there is a functor 
$$\chi: \textbf{sSet}/A \rightarrow [A, \textbf{sSet}^\sharp]$$
from the slice category over $A$. 

\begin{theorem} \label{class_simp}
Let $A$ be a simplicial set. The pair of functors $(\chi, \textbf{dcolim})$ 
$$\chi: \textbf{sSet}/A \leftrightarrows [A, \textbf{sSet}^\sharp] : \textbf{dcolim}$$
is an equivalence of categories.
\end{theorem}

\begin{redi_proof}
Let $f: X \rightarrow A$ be an object in $\textbf{sSet}/A$. By Proposition \ref{smodel} we have the following
\begin{align*}
    \textbf{dcolim}(\chi_f)_n &\cong 
    \{(\sigma, x) | \sigma \in A_n, x \in \chi_f(\sigma)_n \ \text{section} \} \\
    &= \{(\sigma, x) | \sigma \in A_n, x \in X_n \text{and} f(x) = \sigma  \} \\
    &= \{ (f(x), x) | x \in X_n \} \\
    &\cong X_n 
\end{align*}
It is easy to see that these isomorphisms of sets assemble to an isomorphism $\textbf{dcolim}(\chi_f) \cong X$.

On the other hand, let $F : A \rightarrow \textbf{sSet}^\sharp$ be a horizontal diagram. Let $X$ be the simplicial set defined in Construction \ref{X} and $f : X \rightarrow A$ be the map defined in Corollary \ref{slice}. $(X,f)$ is a model for the double colimit of $F$ by Proposition \ref{smodel}. We prove that $\chi_f \cong F$ as horizontal diagrams.

Let $\sigma \in A_n$ be an $n$-simplex of $A$. By definition $\chi_f(\sigma) \cong X_\sigma$ which is given by the pullback

\begin{center}
    \begin{tikzpicture}

\node (v1) at (-2.5,3) {$X_\sigma$};
\node (v4) at (-2.5,1) {$\Delta^n$};
\node (v2) at (-0.5,3) {$X$};
\node (v3) at (-0.5,1) {$A$};
\draw[->, dashed]  (v1) edge (v2);
\draw[->]  (v2) edge node[right, font = \scriptsize]{$f$} (v3);
\draw[->, dashed]  (v1) edge  (v4);
\draw[->]  (v4) edge node[above, font = \scriptsize]{$\sigma$} (v3);
\end{tikzpicture}
\end{center}

Fix $[k] \in \Delta$. Then a $k$-simplex of $X_\sigma$ is represented by a pair $(\tau, (\eta, x))$ where $\tau : [k] \rightarrow [n]$ is a map in $\Delta$ and $(\eta, x)$ is a $k$-simplex of $X$. Recall from Construction \ref{X} that the latter signify a $k$-simplex $\eta \in A_k$ and a section $x \in F(\eta)_k$ of $p_\eta$. These must satisfy $\sigma \theta = f(\eta, x)$ or simply $\sigma \theta = \eta$ by definition of $f$. Thus we conclude $$X_{\sigma, k} \cong \{(\tau, x) \ | \tau: [k] \rightarrow [n] \ \text{in} \ \Delta, x \in F(\sigma \tau)_k \ \text{section of} \ p_{\sigma \theta} \}$$

On the other hand, we have a pullback square

\begin{center}
    \begin{tikzpicture}

\node (v1) at (-2.5,3) {$F(\sigma \tau)$};
\node (v4) at (-2.5,1) {$\Delta^k$};
\node (v2) at (-0.5,3) {$F(\sigma)$};
\node (v3) at (-0.5,1) {$\Delta^n$};
\draw[->, dashed]  (v1) edge (v2);
\draw[->]  (v2) edge node[right, font = \scriptsize]{$p_\sigma$} (v3);
\draw[->, dashed]  (v1) edge  (v4);
\draw[->]  (v4) edge node[above, font = \scriptsize]{$\tau$} (v3);
\end{tikzpicture}
\end{center}

Hence the data for the section $x \in F(\sigma \tau)_k$ is simply the data of a $k$-simplex $y : \Delta^k \rightarrow F(\sigma)$ with $\tau = p_\sigma y$. But this data is determined solely by the choice of the $k$-simplex $y$
and therefore we conclude 
$$X_{\sigma, k} \cong F(\sigma)_k$$
for all $\sigma$. Again, it is easy to verify that these isomorphisms assemble to define an isomorphism of diagrams $\chi_f \cong F$
\begin{flushright}
$\square$
\end{flushright}
\end{redi_proof}

\begin{corollary} \label{class_inner}
Let $A$ be a simplicial set. The pair of functors $(\chi, \textbf{dcolim})$ 
$$\chi: (\textbf{sSet}/A)_{inner} \leftrightarrows [A, \infty\textbf{Cat}^\sharp] : \textbf{dcolim}$$
is an equivalence of categories between the full subcategory of $\textbf{sSet}/A$ consisting of inner fibrations and the category of $A$-indexed diagrams of correspondences of $\infty$-categories. Here double colimits are computed in $\textbf{sSet}^\sharp$.
\end{corollary}

\begin{redi_proof}
A simplicial map $p: X \rightarrow A$ is an inner fibration if and only if $X_\sigma$ is an $\infty$-category for all $\sigma \in A_n$ and $[n] \in \Delta$, i.e. if the classifying diagram $\chi_p : A \rightarrow \textbf{sSet}^\sharp$ from Theorem \ref{class_simp} factors through the inclusion $\infty\textbf{Cat}^\sharp \subseteq \textbf{sSet}^\sharp$. 
\begin{flushright}
$\square$
\end{flushright}
\end{redi_proof}

\begin{note}
It is reasonable to expect a higher version of Theorem \ref{class_simp} to hold. That would be a theorem which takes into consideration the fact that $\textbf{sSet}/A$ is simplicially enriched for all simplicial sets $A$. We briefly sketch how this can be achieved.

There is a natural simplicial enrichment of the diagram category $[A, \textbf{sSet}^\sharp]$ as well. Given a diagram $F\in [A, \textbf{sSet}^\sharp]$ define $F \times \Delta^n$ to be the composite
$$F \times \Delta^n : A \times \Delta^n \rightarrow A \xrightarrow{F} \textbf{sSet}^\sharp$$
of $F$ with the projection map. Then, given $F,G \in [A, \textbf{sSet}^\sharp]$ we may define a simplicial set $[A, \textbf{sSet}^\sharp](F,G)$ with set of $n$-simplices given by the set of vertical transformations $F \times \Delta^n \Rightarrow G \times \Delta^n$
$$[A, \textbf{sSet}^\sharp](F,G)_n = \textbf{vTrans}(F \times \Delta^n, G \times \Delta^n)$$

It should not be difficult, perhaps using tools which are a little bit more abstract than what we used here, to show that there is an equivalence of \textbf{sSet}-categories
$$\textbf{sSet}/A \cong [A, \textbf{sSet}^\sharp]$$
Hopefully, we will give a detailed proof in future work.
\end{note}

\begin{remark}

In particular, we believe it would be interesting to define the ``double $\infty$-category" of $\infty$-categories mentioned in \cite{barwick2018fibrations} as a certain bisimplicial set and carry out a theory of double colimits in a way which is analogous to the theory of colimits in $\infty$-category theory.
\end{remark}



\end{document}